\documentclass[aos,MSNbibl,seceqn,nameyear,dvips]{arximspdf}
\usepackage{graphicx}
\doi{10.1214/13-AOS1168} %kopijuoti is PTS
\volume{41}
\issue{6}
\pubyear{2013}
\firstpage{2820}
\lastpage{2851}

\makeatletter

\newtheorem{theorem}{Theorem}
\newtheorem{proposition}{Proposition}

\usepackage{mathrsfs}
\makeatother

\begin{document}
\begin{frontmatter}

\title{Tests alternative to higher criticism for high-dimensional means
under sparsity and column-wise dependence}
\runtitle{Tests alternative to higher criticism}

\begin{aug}
\author[a]{\fnms{Ping-Shou} \snm{Zhong}\ead[label=e1]{pszhong@stt.msu.edu}\thanksref{tt1}},
\author[b]{\fnms{Song Xi} \snm{Chen}\corref{}\ead[label=e2]{songchen@iastate.edu}\thanksref{tt2,tt4}}
\and
\author[c]{\fnms{Minya} \snm{Xu}\ead[label=e3]{minyaxu@gsm.pku.edu.cn}\thanksref{tt3,tt4}}
\thankstext{tt1}{Supported in part by NSF Grants DMS-13-09156 and DMS-12-09112.}
\thankstext{tt2}{Supported in part by NSF Grant DMS-13-09210 and NSFC Key Grant
10901010 and Grant 71371016.}
\thankstext{tt3}{Supported in part by NSFC Key Grant 10901010 and the Center for
Statistical Science at Peking University.}
\thankstext{tt4}{Supported in part by Key Laboratory of Mathematical Economics and
Quantity Finance of MOE at Peking University.}
\runauthor{P.-S. Zhong, S. X. Chen and M. Xu}
\affiliation{Michigan State University, Peking University and Iowa
State University, and Peking University}
\address[a]{P.-S. Zhong\\
Department of Statistics and Probability\\
Michigan State University\\
East Lansing, Michigan 48824\\
USA\\
\printead{e1}}

\address[b]{S. X. Chen\\
Department of Business Statistics\\
\quad and Econometrics\\
Guanghua School of Management\\
\quad and Center for Statistical Science\\
Peking University\\
Beijing 100871\\
China\\
and\\
Department of Statistics\\
Iowa State University\\
Ames, Iowa 50011-1210 \\
USA\\
\printead{e2}}

\address[c]{M. Xu\\
Department of Business Statistics\\
\quad and Econometrics\\
Guanghua School of Management\\
\quad and Center for Statistical Science\\
Peking University\\
Beijing 100871\\
China\\
\printead{e3}}
\end{aug}

% HISTORY:
\received{\smonth{3} \syear{2013}}
\revised{\smonth{9} \syear{2013}}

% ABSTRACT
%
\begin{abstract}
We consider two alternative tests to the Higher Criticism test of
Donoho and Jin [\textit{Ann. Statist.} \textbf{32} (2004) 962--994] for
high-dimensional means under the sparsity of the
nonzero means for sub-Gaussian distributed data with unknown
column-wise dependence.
The two alternative test statistics are constructed by first
thresholding $L_1$ and $L_2$ statistics based on the sample means, respectively,
followed by maximizing over a range of thresholding levels to make the
tests adaptive to the unknown signal strength and sparsity. The two
alternative tests can attain
the same detection boundary of the Higher Criticism test in
[\textit{Ann. Statist.} \textbf{32} (2004) 962--994]
which was established for uncorrelated Gaussian data.
It is demonstrated that the maximal $L_2$-thresholding test is at least
as powerful as the maximal $L_1$-thresholding test, and both the
maximal $L_2$ and $L_1$-thresholding tests are at least as powerful as
the Higher Criticism test.
\end{abstract}

% KEYWORDS
\begin{keyword}[class=AMS]
\kwd[Primary ]{62H15}
\kwd[; secondary ]{62G20}
\kwd{62G32}
\end{keyword}
\begin{keyword}
\kwd{Large deviation}
\kwd{large $p$}
\kwd{small $n$}
\kwd{optimal detection boundary}
\kwd{sparse signal}
\kwd{thresholding}
\kwd{weak dependence}
\end{keyword}

\end{frontmatter}

%s1 #&#
\section{Introduction}\label{sec1}

Let $\mathbf{X}_1,\ldots,\mathbf{X}_n$ be independent and identically
distrib\-uted (I.I.D.) $p$-variate
random vectors generated from the following model:
%
%e1.1 #&#
\begin{equation}
\label{model} \mathbf{X}_i=\mathbf{W}_i+\bolds{\mu}\qquad
\mbox{for } i=1,\ldots ,n,
\end{equation}
where
$\bolds{\mu}=(\mu_{1},\ldots,\mu_{p})^T$ is a $p$-dimensional unknown
vector of means, $\mathbf{W}_i=(W_{i1},\ldots, W_{i p})^T$ and $\{
\mathbf{W}_i\}_{i=1}^n$ are I.I.D. random vectors with zero mean and
common covariance $\bolds{\Sigma}$. For the $i$th sample, $\{W_{ij}\}
_{j=1}^p$ is a sequence of weakly stationary dependent random variables
with zero mean and
variances $\sigma_j^2$. Motivated by the high-dimensional applications
arising in genetics, finance and other fields,
the current paper focuses on testing high-dimensional hypotheses
%
%e1.2 #&#
\begin{equation}
\label{nullhyper} H_0\dvtx \bolds{\mu}=0\quad \mbox{vs} \quad H_1\dvtx
\mbox{nonzero $\mu_j$ are sparse and faint.} %not 0 for some $
\end{equation}
The specifications for the sparsity and faintness in the above $H_1$
are the following. There are $p^{1-\beta}$ nonzero $\mu_j$'s (signals)
for a $\beta\in(1/2,1)$, which are sparse since %are sparsely
%populated such that
the signal bearing dimensions constitute only a small fraction of the
total $p$ dimensions. Also under the $H_1$, the signal strength is
faint in that the nonzero $\mu_j = \sqrt{ 2 r \log(p)/n}$ for $r \in
(0, 1)$.
%strength through the rest of this paper.}
These specification of the $H_1$ have been the most challenging
``laboratory'' conditions in developing novel testing procedures under
high dimensionality.

\citet{DonohoJin} pioneered the theory of the Higher Criticism (HC)
test which was originally conjectured in \citet{Tukey}, and showed that
the HC test can attain the optimal detection boundary established by
\citet{Ingster} for uncorrelated Gaussian random vectors ($\bolds
{\Sigma
}=\mathbf{I}_p$).
The optimal detection boundary is a phase-diagram in the space of
$(\beta, r)$, the two quantities which define the sparsity and the
strength of nonzero $\mu_j$'s under the $H_1$,
such that if $(\beta, r)$ lies above the boundary, there exists a test
which has asymptotically diminishing probabilities of the type I and
type II errors simultaneously; and if $(\beta, r)$ is below the
boundary, no such test exists.
%The investigation of \citet{DonohoJin} was conducted for an
%independent and identically distributed (IID) random vector, namely,
%for a %column-wise independence vector ($n=1$).
Hall and Jin (\citeyear{HallJin}, \citeyear{HallJin2010}) investigated the impacts of the column-wise
dependence on the HC test. In particular, \citet{HallJin} found
that the HC test is adversely affected if the dependence is of long
range dependent. If the dependence is weak, and the covariance matrix
is known or can be estimated reliably, the dependence can be utilized
to enhance the signal strength of the testing problem so as to improve
the performance of the HC test. The improvement is reflected in
lowering the needed signal strength $r$ by a constant factor. \citet{Delaigle2009} evaluated the HC test under a nonparametric setting
allowing column-wise dependence, and showed that {the detection
boundary of \citet{DonohoJin}} for the HC test can be maintained
under weak column-wise dependence. \citet{Delaigle} showed
that the standard HC test based on the normality assumption can perform
poorly when the underlying data deviate from the normal distribution
and studied a version of the HC test based on the $t$-statistics
formulation. \citet{Cai2011} considered detecting Gaussian
mixtures which differ from the null in both the mean and the variance.
%%for detecting heterogeneous mixtures,
Arias-Castro, Bubeck and Lugosi (\citeyear{Arias-Castro2012a,Arias-Castro2012b}) established the lower
and upper bounds for the minimax risk for detecting sparse differences
in the covariance.

We show in this paper that there are alternative test procedures for
weakly dependent sub-Gaussian data with unknown covariance which attain
the same detection boundary as the HC test established in \citet{DonohoJin} for Gaussian distributed data with $\bolds{\Sigma
}=\mathbf{I}_p$.
The alternative test statistics are obtained by first constructing,
for $\gamma=1$ and $2$,
\[
T_{\gamma n}(s) = \sum_{j=1}^p |
\sqrt{n}\bar{X}_j/\sigma_j|^{\gamma
} I\bigl( |
\bar{X}_j| \geq\sigma_j\sqrt{\lambda_p(s)/n}
\bigr),
\]
which threshold with respect to $\bar{X}_j$ at a level $\sqrt{\lambda
_p(s)/n}$ for $s \in(0,1)$, where $\lambda_p(s)=2s\log p$, $\bar
{X}_j$ is the sample mean of the $j$th margin of the data vectors and
$I(\cdot)$ is the indicator function. We note that $\gamma=1$ and $2$
correspond to the $L_1$ and $L_2$ versions of the thresholding
statistics, respectively; and $\gamma=0$ corresponds to the HC test
statistic. In the literature, the $L_1$ statistic is called the hard
thresholding in \citet{DonohoJohnstone} and \citet{DonohoJin08},
and the $L_0$ statistic is called the clipping thresholding in \citet{DonohoJin08}.
%Prof. Chen, They referred to nonmaximal version. They discussed the
%choice of the %thresholds in multiple testing using HC with clipping,
%hard and soft versions. } }
%, and $\gamma=1$ and $2$ representing respectively the $L_1$ and
%$L_2$ versions.
We then %standardize $T_{\gamma n}(s)$
maximize standardized versions of $T_{\gamma n}(s)$ with respect to $s$
over $\mathcal{S}$, a subset of $(0,1)$, which results in the following
maximal $L_\gamma$-thresholding statistics:
%
%e1.3 #&#
\begin{equation}
\hat{\mathcal{M}}_{\gamma n}=\max_{s\in\mathcal{S}}
\frac
{T_{\gamma
n}(s)-\hat{\mu}_{T_{\gamma n},0}(s)}{\hat{\sigma}_{T_{\gamma n},0}(s)}\qquad \mbox{for $\gamma=0, 1$ and $2$,} \label{eq:Mgamman}
\end{equation}
where %$\mathcal{S}$ is a subset of $(0,1)$, and
$\hat{\mu}_{T_{\gamma n},0}(s)$ and $\hat{\sigma}_{T_{\gamma n},0}(s)$
are, respectively, estimators of the mean\break ${\mu}_{T_{\gamma n},0}(s)$
and standard deviation ${\sigma}_{T_{\gamma n},0}(s)$ of $T_{\gamma
n}(s)$ under $H_0$, whose forms will be given later in the paper.
%We will evaluate the test performance for two values of $\gamma$: $
%We would like to explore the opportunity offered by other values of $
{By developing the asymptotic distributions of $\hat{\mathcal
{M}}_{\gamma n}$, the maximal $L_{\gamma}$-thresholding tests are
formulated for $\gamma=0, 1$ and $2$ with the maximal $L_0$-test being
equivalent to the HC test.} An analysis on the relative power
performance of the three tests reveals that if the signal strength
parameter $r \in(0,1)$, the maximal $L_2$-thresholding test is at
least as powerful as the maximal $L_1$-thresholding test, and both the
$L_1$ and $L_2$-thresholding tests are at least as powerful as the HC
test. If we allow a slightly stronger signal so that $r > 2\beta-1$,
the differential power performance of the three tests is amplified with
the maximal $L_2$-test being the most advantageous followed by the
maximal $L_1$-test.

In addition to the connection to the HC test, the maximal $L_{\gamma
}$-thresholding test, by its nature of formulation, is related to the
high-dimensional multivariate testing procedures, for instance, the
tests proposed by \citet{BaiSara} and \citet{ChenQin}.
While these tests can maintain accurate size approximation under a
diverse range of dimensionality and column-wise dependence, their
performance is hampered when the nonzero means are sparse and faint.
The proposed test formulation is also motivated by a set of earlier
works including \citet{DonohoJohnstone} for selecting significant
wavelet coefficients, and \citet{Fan1996} who considered testing for the
mean of a random vector $\mathbf{X}$ with I.I.D. normally distributed components.
We note that the second step of maximization with respect to $s \in
\mathcal{S} \subset(0,1)$ is designed to make the test adaptive to the
underlying signals strength and sparsity, which is the essence of the
HC procedure in \citet{DonohoJin}, as well as that of \citet{Fan1996}.

The rest of the paper is organized as follows. In Section~\ref{sec2} we provide
basic results on the $L_2$-thresholding statistic via the large
deviation method and the asymptotic distribution of the single
threshold statistic. Section~\ref{sec3} gives the asymptotic distribution of
$\hat{\mathcal{M}}_{2n}$ as well as the associated test procedure.
Power comparisons among the HC and the maximal $L_1$ and
$L_2$-thresholding tests are made in Section~\ref{sec4}. Section~\ref{sec5} reports
simulation results which confirm the theoretical results. Some
discussions are given in Section~\ref{sec6}. All technical details are relegated
to the \hyperref[app]{Appendix}.

%s2 #&#
\section{Single threshold test statistic}\label{sec2}

Let $\mathbf{X}_1, \ldots, \mathbf{X}_n$ be an independent
\mbox{$p$-dimensional} random sample from a common distribution $F$,
% where the $i$-th data vector $X_{i} = (X_{1 i}, \cdots, X_{p i})^T$.
and $\mathbf{X}_i = \mathbf{W}_i + \bolds{\mu}$, where
$\bolds{\mu} = (\mu_1, \ldots, \mu_p)^T$ is the vector of means and
$\mathbf
{W}_{i}= (W_{i 1}, \ldots,\break  W_{i p})^{T}$ is a vector consisting of
{potentially dependent random variables} with zero mean and finite
variances. The dependence among $\{W_{i j}\}_{j=1}^p$ is called the
column-wise dependence in $\mathbf{W}_i$. Those nonzero $\mu_j$ are
called ``signals.''

%To make the testing problem the most challenging,
%As already outlined in the description of the $H_1$, we adopt the
%sparse and faint signal setting
%where there are only $m=p^{1-\beta}$ signals for $\beta\in(1/2, 1)$
%and each signal $\mu_j = \sqrt{2 r_j \log(p)/n}$.
%Here $r_i > 0$ represent the signal strength. %, whose range is $r_i
%Although the test statistic $T_{\gamma n}(s)$ can %be used to detect
%both positive and negative signals, we will assume that $r_i>0$ for
%simplifying notation. Indeed, for those dimensions with negative
%signal strength, we could %replace $\bar{X}_j$ by $-\bar{X}_j$ in $T_{
% As \citet{HallJin2010}, %\textbf{we assume the signals' locations %$
% $\ell_1<\ell_2<\cdots<\ell_m$ for $m=p^{1-\beta}$
% are randomly selected from $\{1, 2, \cdots, p\}$ without replacement
%$\mathscr{L}_m=\{(p_1,\cdots,p_m): 1\leq p_1< p_2<\cdots< p_m\leq p
% so that
%P\{\ell_1=p_1,\cdots,\ell_m=p_m\}={p\choose m}^{-1} \mbox{for all %$
%$1 \le p_1 < p_2 < \cdots< n_m \le p$}

Let $\bar{X}_j={n}^{-1}\sum_{i=1}^nX_{i j}$, $\sigma_j^2 =\operatorname
{Var}(W_{i j})$ and $s_{j}^2=(n-1)^{-1}\sum_{i=1}^n(X_{ij}-\bar
{X}_j)^2$ be the sample variance for the $j$th margin.
%are still able to use $s_{j}^2$ to estimate $\sigma_j^2$. Because,
%under the %model (1.1),%$$
%s_j^2=\frac{1}{n-1}\sum_{i=1}^n(X_{ij}-\bar{X}_j)^2=\frac{1}{n-1}
%Only the locations of the signals are random, but the places where
%signal appears are the same for $n$ sample vector. So they were
%cancelled out.
%}}
The signal strength in the $j$th margin can be measured by the
$t$-statistics $\sqrt{n} \bar{X}_j /s_j$ or the $z$-statistics $\sqrt {n} \bar{X}_j /\sigma_j$ if $\sigma_j$ is known.
% Let $Y_{j,n}$ be the square of either the $t$- or z-statistic.
For easy expedition, the test statistics will be constructed based on
the $z$-statistics by assuming $\sigma_{j}$ is known and, without loss of
generality, we assume $\sigma_j^2=1$.
Using the $t$-statistics actually leads to less restrictive conditions
for the underlying random variables since the large deviation results
for the self-normalized
$t$-statistics can be established under weaker conditions to allow
heavier tails in the underlying distribution as demonstrated in \citet{Shao}, \citet{Jing} and \citet{WangHall}.
{See \citet{Delaigle} for analysis on the sparse signal
detection using the $t$-statistics.}

We assume the following assumptions in our analysis:
\begin{longlist}[(C.1)]
\item[(C.1)] The dimension $p = p(n) \to\infty$ as $n \to\infty$ and $\log
(p)=o(n^{1/3})$.

\item[(C.2)] {There exists a positive constant $H$ such that, for any $j \neq
l\in\{1,\ldots, p\}$, }
%the bivariate joint moment generating function of $(W_{j 1}, W_{l 1})$
$E(e^{h^{T}(W_{1 j}^d,W_{1 l}^d)})<\infty$ for $h \in[-H,H]\times[-H,H]$
and $d=2$.

\item[(C.3)] For each $i=1, \ldots,n$, $\{W_{ij}\}_{j=1}^p$ is a {weakly
stationary} sequence such that $E(W_{ij})=E(W_{i(j+k)}) =0$ and $\operatorname
{Cov}(W_{ij},W_{i(j+k)})$ does not depend on $j$ for any integer $k$.
And $\sum_k|\rho_k|<\infty$ where $\rho_k=\operatorname{Cov}(W_{i1}, W_{i(k+1)})$.

%To gain more insight on this aspect, we confine to the following
%situation of
% sparse and faint signal situation: % as specified by the following
%condition:

\item[(C.4)] Among the $p$ marginal means, there are $m=p^{1-\beta}$ signals
for a $\beta\in(1/2, 1)$ and the signal $\mu_j = \sqrt{2r\log(p)/n}$
for a $r > 0$. The signals' locations $\ell_1<\ell_2<\cdots<\ell_m$ are
randomly selected from $\{1, 2, \ldots, p\}$ without replacement
%$\mathscr{L}_m=\{(p_1,\cdots,p_m): 1\leq p_1< p_2<\cdots< p_m\leq p
so that
%
%e2.1 #&#
\begin{eqnarray}
\label{random-location} P(\ell_1=p_1,\ldots,\ell_m=p_m)
=\pmatrix{p\cr  m}^{-1}
\nonumber
\\[-8pt]
\\[-8pt]
\eqntext{\mbox{for all %$
$1 \le
p_1 < p_2 < \cdots< p_m \le p$}.}
\end{eqnarray}
\end{longlist}

(C.1) specifies the growth rate of $p$ relative to the sample size $n$ is
in the paradigm of ``large $p$, small $n$.'' {That $\log p=o(n^{1/3})$
is the rate we can attain for Gaussian data or cases where we can
attain ``accurate'' enough estimation of $\mu_{T_{\gamma n},0}$, which
satisfies equation (\ref{eq:cri1}). When data are not Gaussian and the
``accurate'' estimators are not attainable, the growth rate of $p$ will
be more restrictive at $p=n^{1/\theta}$ ($\theta>0$), as will be
discussed in the next section.}
(C.2) assumes the joint distributions of $(W_{i j}, W_{i l})$ is
sub-Gaussian, which implies each marginal $W_{ij}$ is sub-Gaussian as
well. (C.3) prescribes weak dependence among $\{W_{ij}\}_{j=1}^p$. The
first part of (C.4) reiterates the sparse and faint signal setting. The
range of the signal strength includes the case of $r \in(0,1)$,
representing the most fainted detectable signal strength, which has
been considered in \citet{DonohoJin} and other research works.
%The range also covers $r \ge1$, representing stronger signals.
%Considering $r \ge1$ is to gain a wider view on the relative
%performance of the tests.
The second part of (C.4) provides a random allocation mechanism for the
signal bearing dimensions, which is the same as the one assumed in \citet{HallJin2010}.
Existing research on the detection boundary of the
HC test for the sparse mean problem [\citet{DonohoJin}; \citet{HallJin2010}]
is largely conducted for the case of $n=1$ when the data
are Gaussian. This is understandable since the sample means are
sufficient statistics and there is no loss of {generality} when we
treat the problem as $n=1$, even if we have multiple observations.
However, when the underlying distributions are as specified in (C.2), we
cannot translate the test problem to $n=1$ without incurring a loss of
information.

We first consider the $L_2$ version of the thresholding statistic $T_{2
n}$ %and $\mathcal{M}_{\gamma n}$
in this section. The study of the $T_{1 n}$ version is outlined in
Section~\ref{sec4} when we compare the power performance to the HC test. Let
$Y_{j,n} = n \bar{X}_j^2$. Then,
the $L_2$-thresholding statistic can be written as
%
%e2.2 #&#
\begin{equation}
T_{2 n}(s)=\sum_{j=1}^p
Y_{j,n}I\bigl\{Y_{j,n}\geq\lambda_p(s) \bigr\},
\label{eq:Tns}
\end{equation}
where $s$ is the thresholding parameter that takes values over a range
within $(0,1)$.
%To appreciate this choice of $s$,
%we note that
% from Theorem 5.23 of \citet{Petrov} that the sub-Gaussian assumption
%in C.2 means that, uniformly for $z=o(n^{1/6}),$
%P(\sqrt{n}(\bar{X_i}-\mu_i)\geq z)=\bar{\Phi}(z)\{1+o(1)\}
%where $\bar{\Phi}(\cdot)$ is the standard normal survival function.
%Hence,
% under $H_0$
%P(\max_{1\leq j\leq p} \sqrt{n} \bar{X}_j/\sigma_j \geq\sqrt{2\log
%p})&\leq\sum_{i=1}^p P(\sqrt{n} \bar{X}_j/\sigma_j \geq\sqrt{2\log
%%p})=\sum_{i=1}^p\bar{\Phi}(\sqrt{2\log p})\{1+o(1)\}\\
%&=(2\pi)^{-1/2}(2\log p)^{-1/2}\{1+o(1)\}\to0.
%Therefore,
There is no need to consider $s \ge1$ in the thresholding {since
large deviation results given in \citet{Petrov} imply that under $H_0$$,
P( \max_{1\leq j\leq p}Y_{j,n}\le\lambda_p(s)) \to1.
%&\leq\sum_{i=1}^p P(\sqrt{n} \bar{X}_j \geq\sqrt{2\log p})=
%&=(2\pi)^{-1/2}(2\log p)^{-1/2}\{1+o(1)\}\to0.
$}

Define a set of slowing varying functions:
$L_p^{(1)}=2r\log p+1$, $L_p^{(2)}=\break 2\sqrt{s\log p/\pi}$,
$L_p^{(3)}=s(\sqrt{s}-\sqrt{r})^{-1}\sqrt{\log p/\pi}$, $L_p^{(4)}=
8r\log p$,
%$L_p^{(4)}= 8r\Phi^*_{r,s}\log p$,
$L_p^{(5)}=4s^{3/2}\times \pi^{-{{1}/{2}}}(\log p)^{3/2}$ and
$L_p^{(6)}={2s^2(\log p)^{3/2}}/{\sqrt{\pi}(\sqrt{s}-\sqrt{r})}$.
Let $\phi(\cdot)$ and $\bar{\Phi}(\cdot)$
be the density and survival functions of the standard normal
distribution.

Let $\mu_{T_{2n}, 0}(s)$ %= E\{T_{2n}(s)\}$
and $\sigma^2_{T_{2n}, 0}(s)$ %=Var\{T_{2n}(s)\}$
be the mean and variance of $T_{2n}(s)$ under $H_0$, respectively, and
$\mu_{T_{2n}, 1}(s)$ and
$\sigma^2_{T_{2n}, 1}(s)$ be those, respectively, under the $H_1$ {as
specified in (C.4)}.
The following proposition depicts the mean and variance of $T_{2n}(s)$
{by applying Fubini's theorem and the large deviation results
[\citet{Petrov} and Lemma~A.1 in Zhong, Chen and Xu (\citeyear{ZhenChenXu})].

%pr1 #&#
\begin{proposition}
\label{chap4-cor2} Under \textup{(C.1)--(C.4)}, $E\{T_{2n}(s)\}$ and $\operatorname
{Var}\{T_{2n}(s)\}$ are, respectively,
%
%e2.3 #&#
%e2.4 #&#
\begin{eqnarray}\label{eq:meanTn0}\quad
&&\mu_{T_{2n}, 0}(s)
\nonumber
\\[-8pt]
\\[-8pt]
\nonumber
&&\qquad= p\bigl\{2{\lambda_p^{1/2}(s) }\phi
\bigl({\lambda ^{1/2}_p(s)}\bigr)+2\bar{\Phi}\bigl({
\lambda^{1/2}_p(s)}\bigr)\bigr\} \bigl\{1+O\bigl\{
n^{-1/2}{\lambda^{3/2}_p(s)}\bigr\}\bigr\},
\\
% &=L_p^{(2)}p^{1-s}\{1+o(1)\}; \nn\\
\label{eq:varTn0}&&\sigma_{T_{2n},0}^2 (s)
\nonumber
\\[-8pt]
\\[-8pt]
\nonumber
&&\qquad= p\bigl\{2\bigl[
\lambda^{3/2}_p(s)+3{\lambda ^{1/2}_p(s)}
\bigr]\phi\bigl({\lambda^{1/2}_p(s)}\bigr)+6\bar{\Phi}\bigl({
\lambda ^{1/2}_p(s)}\bigr)\bigr\} \bigl\{1+o(1)\bigr\}
 % &=L_p^{(5)}p^{(1-s)}\{1+o(1)\};\nn\\
\end{eqnarray}
under the $H_0$; and
%combined two terms with same order %$p^{1-\beta-(\sqrt{s}-
\begin{eqnarray*}
%L_p^{(2)}p^{1-s}+L_p^{(3)}p^{1-\beta-(\sqrt{s}-\sqrt{r})^2}I(s>r)\}
\mu_{T_{2n}, 1}(s)&=& \bigl
\{L_p^{(1)}p^{1-\beta}I(s<r)+L_p^{(3)}p^{1-\beta
-(\sqrt{s}-\sqrt{r})^2}I(s>r)
\bigr\} \bigl\{1+o(1)\bigr\}\\
&&{}+\mu_{T_{2n},
0}(s),
\\
\sigma_{T_{2n},1}^2 (s)&=&\bigl\{L_p^{(4)}p^{1-\beta}I(s<r)+
L_p^{(5)}p^{1-s}+L_p^{(6)}p^{1-\beta-(\sqrt{s}-\sqrt{r})^2}I(s>r)
\bigr\}\\
&&{}\times \bigl\{ 1+o(1)\bigr\}
\end{eqnarray*}
under the $H_1$ specified in \textup{(C.4)}.
\end{proposition}

Expressions (\ref{eq:meanTn0}) and (\ref{eq:varTn0}) provide the first
and the second order terms of $\mu_{T_{2n}, 0}(s)$ and $\sigma
_{T_{2n},0}^2 (s)$, which are needed when we consider their empirical
estimation under $H_0$ when formulating the $L_2$ thresholding test
statistic. Note that $\mu_{T_{2n}, 0}(s)=L_p^{(2)}p^{1-s}\{1+o(1)\}$
and $\sigma_{T_{2n},0}^2 (s)=L_p^{(5)}p^{1-s}\{1+o(1)\}$.
Only the first order terms for the variance are needed under $H_1$, but
the approximation to $\mu_{T_{2n},1}(s)$ has to be more accurate so as
to know the order of the difference between $\mu_{T_{2n},1}(s)$ and
$\mu
_{T_{2n},0}(s)$.
%The order of the relative approximation error using the leading order
%term is also given in (\ref{eq:meanTn0}).}
Proposition \ref{chap4-cor2} indicates that the column-wise dependence as specified in
(C.3) does not have much leading order impact on the variance of $T_{2
n}(s)$. The leading order variance is almost the same when $\mathbf
{W}_i$ are column-wise independent. The difference only appears in the
coefficients of the slow-varying functions $L_p^{(4)}$, $L_p^{(5)}$ and
$L_p^{(6)}$,
%under independent case ?: \textbf{The slow varying functions
%$L_p^{(4)}$,
%%$L_p^{(5)}$ and $L_p^{(6)}$ given here are the same as the the slow
%varying functions for independent case, we have ignored some
%negligible part for %the weakly dependent case.} }
while their orders of magnitude remain unchanged.
The reason behind this phenomena is the thresholding.
It can be understood by an analogue for multivariate Gaussian
distributions with nonzero correlation. Despite the dependence in the
Gaussian distribution, exceedances beyond high thresholds are
asymptotically independent [\citet{Sibuya} and \citet{Joe}].\vadjust{\goodbreak}

We now study the asymptotic distribution of $T_{2 n}(s)$ to prepare for
the proposal of the maximal $L_2$-thresholding statistic.
%Recall that $Y_{i,n} = n \bar{X}_i^2$.
% $T_{2n}(s)$ over a range of thresholding levels. We show that this
%formulation will lead to a test procedure that %can attain the optimal
%detection boundary.
Write
\[
T_{2n}(s)= \sum_{j=1}^p
Z_{j,n}(s),
\]
where $Z_{j,n}(s):=Y_{j,n}I\{Y_{j,n}> \lambda_p(s)\}$ and $\lambda_p(s)
=2 s \log(p)$.
For integers $a, b \in[-\infty,\infty]$ such that $a < b$, define
$\mathscr{F}_{a}^b=\sigma\{Z_{l,n}(s)\dvtx l\in(a, b)\}$ as the $\sigma
$-algebra generated by $\{Z_{l,n}(s)\}_{l=a}^b$ and define the $\rho
$-mixing coefficients
%
%e2.5 #&#
\begin{equation}
\rho_{Z(s)}(k) = %\rho(\mathscr{F}_{-\infty}^i,\mathscr{F}^{
\sup_{l, \xi\in L^2(\mathscr{F}_{-\infty}^l ), \zeta\in
L^2(\mathscr
{F}^{\infty}_{l+k} ) }\bigl|
\operatorname{Corr}(\xi,\zeta)\bigr|.\label{eq:mixing}
\end{equation}
See \citet{Doukhan} for comprehensive discussions on the mixing concept.
The following is a condition regarding the dependence among $\{
Z_{j,n}(s)\}_{j=1}^p$.

\begin{longlist}[(C.5)]
\item[(C.5)] For any $s \in(0,1)$, the sequence of random variables $\{
Z_{j,n}(s)\}_{j=1}^{p}$ is $\rho$-mixing such that
%the mixing coefficients $\rho_{Z}(k)$ satisfying
$\rho_{Z(s)}(k)\leq C\alpha^k$ for some $\alpha\in(0,1)$ and a
positive constant~$C$. %$\lim_{k\to\infty}\rho_{Z(s)}(k)=0$.
\end{longlist}

The requirement of $\{Z_{j,n}(s)\}_{j=1}^p$ being $\rho$-mixing for
each $s$ is weaker than requiring the original data
columns $\{X_{i j}\}_{j=1}^p$ being $\rho$-mixing, whose mixing
coefficient $\rho_{X_i}(k)$ can be similarly defined as (\ref{eq:mixing}).
%In fact, that $X$ being $\rho$-mixing implies the $\rho$-mixing of $\{
%Z_{j,n}(s) \}_{j=1}^p$ for any $s \in(0,1)$.
This is because, according to Theorem 5.2 in \citet{Bradley},
\[
\rho_{Z(s)}(k)\leq\sup_{i\leq n}\rho_{X_i}(k)=
\rho_{X_1}(k)\qquad \mbox{for each $k=1,\ldots,p$ and $s \in(0,1)$.}
\]
%
%where $\rho_{X}(k)=\rho_{X_j}(k)$ for $j=1,\cdots, n$.

%Denote the mean and variance of $T_{2n}(s)$ under the $H_1$ as $
%%$\sigma^2_{T_{2n} 1}(s)$ respectively, whose leading order terms can
%be deducted from
%Corollary 2.
The following theorem reports the asymptotic normality of $T_{2n}(s)$
under both $H_0$ and $H_1$.

%th1 #&#
\begin{theorem}\label{th1}
Assume \textup{{(C.1)--(C.5)}}.
%need C.3. Because we should be able to approximate the moments of
%$T_{2n}(s)$. We %also need C.5 for the hypothesis considered under
%alternative. I have now added these two conditions.}}
Then, for any $s \in(0, 1)$, %-\eta]$ where $\eta$ is an arbitrary
%constant in $(0,1)$,
\begin{eqnarray*}
&&\hphantom{i}\mathrm{(i)}\quad \sigma^{-1}_{T_{2n}, 0}(s)\bigl\lbrace
T_{2n}(s) - \mu _{T_{2n},0}(s) \bigr\rbrace\stackrel{d} {\to}
N(0,1) \qquad\mbox{under $H_0$};
\\
&&\mathrm{(ii)}\quad \sigma^{-1}_{T_{2n}, 1}(s) \bigl\lbrace
T_{2n}(s) - \mu _{T_{2n}, 1}(s) \bigr\rbrace\stackrel{d} {\to}
N(0,1) \qquad\mbox{under $H_1$}.
\end{eqnarray*}
\end{theorem}

From (\ref{eq:meanTn0}) and (\ref{eq:varTn0}), define the leading order
terms of ${\mu}_{T_{2n}, 0}(s)$ and $\sigma_{T_{2n}, 0}^2(s)$, respectively,
\begin{eqnarray*}
\tilde{\mu}_{T_{2n}, 0}(s)&=& p\bigl\{2{\lambda^{1/2}_p(s)
}\phi\bigl({\lambda ^{1/2}_p(s)}\bigr)+2\bar{\Phi}\bigl({
\lambda^{1/2}_p(s)}\bigr)\bigr\} \qquad\mbox{and}
\\
\tilde{\sigma}_{T_{2n}, 0}^2(s)&=& p\bigl\{2\bigl[
\lambda^{3/2}_p(s)+3{\lambda ^{1/2}_p(s)}
\bigr]\phi\bigl({\lambda^{1/2}_p(s)}\bigr)+6\bar{\Phi}\bigl({
\lambda ^{1/2}_p(s)}\bigr)\bigr\}.
\end{eqnarray*}
It is clear that the asymptotic normality in Theorem \ref{th1}(i) remains if we
replace $\sigma_{T_{2n}, 0}(s)$ by $\tilde{\sigma}_{T_{2n}, 0}(s)$.

To formulate a test procedure based on the thresholding statistic $T_{2
n}(s)$, we need to estimate $\mu_{T_{2n},0}(s)$
by a $\hat{\mu}_{T_{2n}, 0}(s)$, say. Ideally, if
%
%e2.6 #&#
\begin{equation}
\mu_{T_{2n},0}(s) - \hat{\mu}_{T_{2n},0}(s) = o\bigl\{\tilde{\sigma
}_{T_{2n}, 0}(s)\bigr\}, \label{eq:cri1}\vadjust{\goodbreak}
\end{equation}
the first part of Theorem \ref{th1} remains valid if we replace $\mu
_{T_{2n},0}(s)$ with $\hat{\mu}_{T_{2n},0}(s)$. % via the Slutsky
%theorem.
An obvious choice of $\hat{\mu}_{T_{2n},0}(s)$ is $\tilde{\mu}_{T_{2
n}, 0}(s)$, which is known upon given $p$ and $s$.
Indeed, %since $\sqrt{n}\bar{X}_i$ is the standard normal
if $W_{ij}$s are the standard normally distributed,
we have
\[
\mu_{T_{2n},0}(s) = \tilde{\mu}_{T_{2n},0}(s) \qquad\mbox{for } s \in(0,1),
\]
implying the leading order is exactly $\mu_{T_{2n},0}(s)$ for the
Gaussian data. Hence, if we take $\hat{\mu}_{T_{2n},0}(s)=\tilde{\mu
}_{T_{2n},0}(s)$, (\ref{eq:cri1}) is satisfied for the Gaussian data.

For non-Gaussian observations, the difference between $\mu
_{T_{2n},0}(s)$ and\break  $\tilde{\mu}_{T_{2 n}, 0}(s)$ may not be a smaller
order of $\sigma_{T_{2n}, 0}(s)$.
Specifically, from (\ref{eq:meanTn0}) and (\ref{eq:varTn0}), we have
\[
\frac{ \mu_{T_{2n},0}(s) - \tilde{\mu}_{T_{2
n},0}(s)}{\sigma
_{T_{2n},0}(s)}= O \bigl\{{\lambda^{5/4}_p(s)}p^{(1-s)/2}n^{-1/2}
\bigr\}. %=O((
\]
To make the above ratio diminishing to zero, the strategy of \citet{Delaigle} can be adopted by restricting $p=n^{1/\theta}$ and
{ $s \in((1-\theta)_{+}, 1)$ for a positive~$\theta$}, where $(a)_{+}
= a$ if $a > 0$ and $(a)_{+} = 0$ if $a \le0$.
Under this circumstance,
%
%e2.7 #&#
\begin{equation}
\frac{ \mu_{T_{2n},0}(s) - \tilde{\mu}_{T_{2
n},0}(s)}{\sigma
_{T_{2n},0}(s)} =O \bigl\{({2s/\theta\log n})^{5/4}n^{{(1-s-\theta)}/{(2\theta)
}}
\bigr\}\to0. \label{eq:restrict}
\end{equation}
Clearly, for a not so high dimension with $\theta\ge1$, (\ref
{eq:restrict}) holds for all $s \in(0,1)$, % without restricting $s>(1-
and $\tilde{\mu}_{T_{2 n},0}(s)$ satisfies (\ref{eq:cri1}).
For higher dimensions with $\theta< 1$, the thresholding level $s$ has
to be restricted to ensure (\ref{eq:restrict}).
%As will be demonstrated later,
The restriction can
alter the detection boundary of the test we will propose in the next
section. This echoes a similar phenomena for the HC test given in
\citet{Delaigle}.
To expedite our discussion, we assume in the rest of the paper that
(\ref{eq:cri1}) is satisfied by the $\hat{\mu}_{T_{2n},0}(s)$. We note
such an arrangement is not entirely unrealistic, as a separate effort
may be made to produce more accurate estimators. %, for instance by
%utilizing large deviation expansions to the mean.
Assuming so allows us to stay focused on the main agenda of the
testing problem.

The asymptotic normality established in Theorem \ref{th1} allows an asymptotic
\mbox{$\alpha$-level} test
that rejects $H_0$ if
%
%e2.8 #&#
\begin{equation}
T_{2n}(s)-\hat{\mu}_{T_{2n},0}(s)>z_\alpha\tilde{
\sigma}_{T_{2n},
0}(s), \label{eq:test1}
\end{equation}
where $z_\alpha$ is the upper $\alpha$ quantile of the standard normal
distribution.

%s3 #&#
\section{Maximal thresholding}\label{sec3}
While the asymptotic normality of $T_{2n}(s)$ in Theorem \ref{th1} ensures the
single thresholding level test in (\ref{eq:test1}) a correct size
asymptotically, the power of the %above single threshold level
test depends on $s$, the underlying signal strength $r$ and the
sparsity $\beta$.
A test procedure is said to be able to separate a pair of null and
alternative hypotheses % $H_0$ and $H_1$
asymptotically if the sum of the probabilities of the type I and type~II errors converges to zero as $n \to\infty$.
Let $\alpha_n$ be a sequence of the probabilities of type I error,
which can be made converging to zero as $n \to\infty$.
The sum of the probabilities of the type I and type II errors for the test
given in (\ref{eq:test1}) with nominal size $\alpha_n$ is approximately
%
%e3.1 #&#
\begin{equation}
\mathrm{Err}_{\alpha_n}:=\alpha_n+P \biggl(\frac{T_{2n}(s)-\mu
_{T_{2n},0}(s)}{\sigma
_{T_{2n},0}(s)}\leq
z_{\alpha_n} \Big| H_1 \biggr), %:=\alpha_n+Err_{II}.
\label{eq:errors}
\end{equation}
which is attained based on the facts that (i) the size $\alpha_n$ is
attained asymptotically and (ii) $\hat{\mu}_{T_{2n}, 0}(s)$ and
$\tilde
{\sigma}_{T_{2n}, 0}(s)$ are sufficiently accurate estimators { in the
test procedure (\ref{eq:test1}).}

Our strategy is to first make $\alpha_n\to0$ such that $z_{\alpha
_n}=C(\log p)^\varepsilon$ for an arbitrarily
small $\varepsilon>0$ and a constant $C>0$. % and then to consider the
%type II error term,
The second term on the right-hand side of (\ref{eq:errors}) is
%e3.2 #&#
\begin{eqnarray}
&&\mathrm{Err}_{\mathit{II}}:=P \biggl(\frac{T_{2n}(s)-\mu_{T_{2n},1}(s)}{\sigma_{T_{2n},
1}(s)}
\nonumber
\\[-8pt]
\\[-8pt]
\nonumber
&&\hspace*{34pt}\qquad\leq z_{\alpha_n}
\frac{ \sigma_{T_{2n}, 0}(s)}{\sigma_{T_{2n},
1}(s)}-\frac{\mu_{T_{2n},1}(s)-\mu_{T_{2n},0}(s)}{ \sigma_{T_{2n},
1}(s) } \biggr).
\end{eqnarray}
Because $z_{\alpha_n}$ is slowly varying, $0<\sigma
_{T_{2n},0}(s)/\sigma
_{T_{2n},1}(s)\leq1$
and $(T_{2n}(s)- \mu_{T_{2n},1}(s))/\sigma_{T_{2n}, 1}(s)$ is
stochastically bounded, %as implied by Theorem 1, it can be shown that
a necessary and sufficient condition that
ensures $\mathrm{Err}_{\alpha_n}\to0$ is
%
%e3.3 #&#
\begin{equation}
\label{detectable-condition} \Delta_2(s;r,\beta):=\frac{\mu_{T_{2n},1}(s)-\mu
_{T_{2n},0}(s)}{\sigma
_{T_{2n}, 1}(s)}\to\infty.
\end{equation}

From Proposition \ref{chap4-cor2},
it follows that, up to a factor $1+o(1)$,
\begin{eqnarray*}
\Delta_2(s;r,\beta)=\cases{ %
C_1p^{(1+s-2\beta)/2}, %\{1+o(1)\}
& \quad$\mbox{if } s\leq r \mbox{ and } s\leq
\beta$;
\vspace*{2pt}\cr
C_2p^{(1-\beta)/2}, %\{1+o(1)\}
&\quad $\mbox{if } s\leq r \mbox{ and } s>
\beta;$
\vspace*{2pt}\cr
C_3p^{1/2-\beta+r-(\sqrt{s}-2\sqrt{r})^2/2}, %\{1+o(1)\}
& \quad$\mbox{if } s> r \mbox{ and } s\leq(
\sqrt{s}-\sqrt{r})^2+\beta;$
\vspace*{2pt}\cr
C_4p^{(1-\beta-(\sqrt{s}-2\sqrt{r})^2)/2}, %\{1+o(1)\}
&\quad $\mbox{if } s> r \mbox{ and } s> (
\sqrt{s}-\sqrt{r})^2+\beta,$}
\end{eqnarray*}
where $C_1=\sqrt{2}(\pi s)^{{1}/{4}}(\frac{r}{s})(\log
p)^{{1}/{4}}$,
$C_2=\frac{1}{2}(r\log p)^{1/2}$, $C_3={s^{1/4}(\log p)^{-
{1}/{4}}}/\break  \{
\sqrt{2}\pi^{1/4}(\sqrt{s}-\sqrt{r})\}$ and $C_4=(2\sqrt{\pi
}(\sqrt {s}-\sqrt{r}))^{-{{1}/{2}}}(\log p)^{-{1}/{4}}$.

%$C_1=\sqrt{2}(\pi s)^{\fourth}(\frac{r}{s})\Phi_{r,s}^*(\log p)^{

Let
\begin{eqnarray*}
\varrho^{\ast}(\beta)=\cases{%
\beta-1/2, & \quad$\mbox{$1/2 < \beta\leq3/4$;}$
\vspace*{2pt}\cr
(1-\sqrt{1-\beta})^2, &\quad $\mbox{$3/4<\beta<1$.}$}
\end{eqnarray*}
As demonstrated in \citet{DonohoJin} and {\citet{Ingster}},
%%\citet{Delaigle2009},
the phase diagram $r=\varrho^{\ast}(\beta)$ is the optimal detection
boundary for testing the hypotheses we are considering in this paper
when the data are Gaussian and $\bolds{\Sigma}=\mathbf{I}_p$.
Here the optimality means that for any $r > \varrho^{\ast}(\beta)$,
there exists at least one test such that the sum of the probabilities
of the type I and type II errors diminishes to zero as $n \to\infty$;
but for $r < \varrho^{\ast}(\beta)$, no such test exists.
For correlated Gaussian data such that $\bolds{\Sigma} \ne\mathbf
{I}_p$, \citet{HallJin2010} found that the detection boundary
$r=\varrho
^{\ast}(\beta)$ may be lowered by transforming the data via the inverse
of Cholesky factorization $\mathbf{L}$ such that $\mathbf{L}\Sigma
\mathbf{L}^T=\mathbf{I}_p$. {More discussion on the optimality is given
in Section~\ref{sec6}.}

From the expression of $\Delta_2(s;r,\beta)$ given above, it can be
shown (see the proof of Theorem \ref{detect-upper-bound} in the \hyperref[app]{Appendix}) that if $r>\varrho
^{\ast}(\beta)$
there exists at least one $s \in(0,1)$ for each pair of $(r,\beta)$
such that (\ref{detectable-condition}) is satisfied and, hence, the
thresholding test would be powerful. This is the key for the maximal
$L_2$-thresholding test that we will propose later to attain the
detection boundary.

It is clear that we have to make the thresholding level $s$ adaptive to
the unknown $r$ and $\beta$. %As $r$ and $\beta$ are unknown,
One strategy is to use a range of thresholding levels, say, $s \in
{\mathcal{S}} \subset(0,1)$, so that the underlying $(r,\beta)$ can be
``covered.'' This is the very idea of the HC test. % in Donoho and Jin
%(2004).

Let $\hat{\mathcal{T}}_{2,n}(s)= \tilde{\sigma}_{T_{2n},0}^{-1}(s)\{
T_{2n}(s)-\hat{\mu}_{T_{2n},0}(s)\}$ be the standardized version of
$T_{2n}(s)$.
Define the maximal thresholding statistic
%maximizing ${\cal{T}}_2(s)$ over $s \in,
%enough. Can we just choose $S=(0, 1-\eta)$ throughout the paper?: Yes,
%we are %acctually using ${\cal{S}}=(0,1-\eta)$ throughout the paper.}
\[
\hat{\mathcal{M}}_{2 n} = \sup_{s \in{\mathcal{S}}} \hat {\mathcal
{T}}_{2,n}(s), %= \max_{s \in{\cal{S}}} \frac{T_{2n}(s)-\mu_{T_{2n},0}(s)}{
\]
where ${\mathcal{S}}=(0,1-\eta]$ for an arbitrarily small positive
$\eta$.
%which we call the maximal $L_2$-thresholding statistic.
Let
%
%e3.4 #&#
\begin{equation}\qquad
\mathcal{S}_n={\bigl\{s_i\dvtx \mbox{$s_i=Y_{i,n}/(2
\log p)$ and $0<Y_{i,n}<2(1-\eta)\log p$}\bigr\}\cup\{1-\eta
\}}.\label{eq:Sn}
\end{equation}
%
%are between 0 and $2(1-\eta)\log p$, there exists a small interval
%(The last interval) where $\mathcal{T}_n(s)$ attains it maximum at $1-
%restriction.: $1-\eta$ is a point in the set, it is not a restriction
%for $s$. \textbf{Still donot make sense. : $1-\eta$ is an extra point
%in
%addition to those $Y_{i,n}/(2\log p)$. Which simply means that the set
%$S_n$ have one more additional point $1-\eta$. } }
Since both $\hat{\mu}_{T_{2n},0}(s)$ and $\tilde{\sigma}_{T_{2n},0}(s)$
are monotone decreasing functions of $s$, it can be shown that $\hat
{\mathcal{M}}_{2n}$ can be attained on $\mathcal{S}_n$, namely,
%
%e3.5 #&#
\begin{equation}
\label{max-discrete} \hat{\mathcal{M}}_{2n}=\max_{s\in\mathcal{S}_n}
\hat{\mathcal {T}}_{2,n}(s).
\end{equation}
This largely reduces the computational burden of $\hat{\mathcal
{M}}_{2n}$. The asymptotic distribution of $\hat{\mathcal{M}}_{2n}$ is
established in the following theorem.

%th2 #&#
\begin{theorem}
\label{asy-gumbel} Assume \textup{(C.1)--(C.3)}, \textup{(C.5)} and (\ref{eq:cri1})
hold. Then, under~$H_0$,
\[
P\bigl(a(\log p)\hat{\mathcal{M}}_{2 n}-b(\log p,\eta)\leq x\bigr)\to
\exp \bigl(-e^{-x}\bigr), %:=\mathcal{E}
\]
where $a(y)=(2\log(y))^{1/2}$ and $b(y,\eta)=2\log(y)+2^{-1}\log
\log
(y)-2^{-1}\times\break  \log(\frac{4 \pi}{(1-\eta)^2})$.
\end{theorem}

The theorem leads to an asymptotic $\alpha$-level test that rejects
$H_0$ if
%
%e3.6 #&#
\begin{equation}
\hat{\mathcal{M}}_{2 n} >\mathcal{B}_\alpha=\bigl(
\mathcal{E}_\alpha +b(\log p,\eta)\bigr)/a(\log p), \label{eq:L2test}
\end{equation}
where $\mathcal{E}_\alpha$ is the upper $\alpha$ quantile of the Gumbel
distribution $\exp(-e^{-x})$.
We name the test the maximal $L_2$-thresholding test. The following
theorem shows that its detection boundary is $r=\varrho^{\ast}(\beta)$.

%th3 #&#
\begin{theorem}
\label{detect-upper-bound} Under conditions \textup{(C.1)--(C.5)} and
assuming (\ref{eq:cri1}) holds, then
\textup{(i)} if $r>\varrho^{\ast}(\beta)$, %if we reject the null hypothesis
%$H_0$ at some level $\alpha_n\to0$ as $n\to\infty$ such that $z_{
%$C>0$,
the sum of the type I and II errors of the maximal $L_2$-thresholding
tests %with nominal sizes $\alpha_n \to0$
converges to 0 { when the nominal sizes $\alpha_n =\bar{\Phi}((\log
p)^\varepsilon) \to0$ for an arbitrarily} small $\varepsilon>0$ as $n\to
\infty$.
%vague.
%We may still need to specify the rate of $\alpha_n$ converging to 0.
%Because it %does not hold for all $\alpha_n$. }}

\textup{(ii)} If $r<\varrho^{\ast}(\beta)$, the sum of the type I and II errors
of the maximal $L_2$-thresholding test converges to 1 when the nominal
sizes $\alpha_n \to0$ as $n\to\infty$.
\end{theorem}

It is noted that when $r>\varrho^{\ast}(\beta)$ in part (i) of Theorem
\ref{detect-upper-bound}, we need to restrict the rate of the nominal type I error $\alpha
_n$'s convergence to 0, since the conclusion of part (i) may not be
true for all $\alpha_n \to0$. However, in part (ii) where $r <
\varrho
^{\ast}(\beta)$, no restriction for $\alpha_n$ is required, which has
to be the case, as otherwise there is no guarantee that $r=\varrho
^{\ast
}(\beta)$ is the detection boundary of the test.

If the estimator $\hat{\mu}_{T_{2n,0}}(s)$ cannot attain (\ref
{eq:cri1}) and
$\tilde{\mu}_{T_{2n,0}}(s)$ is used as the estimator, we have to
restrict $p=n^{1/\theta}$ for a $\theta\in(0, 1)$ and
limit $s \in(1-\theta, 1)$. In this case, the above theorem is valid
if we replace
$\varrho^{\ast}(\beta)$ by $\varrho^{\ast}_{\theta}(\beta)$, where
\begin{eqnarray*}
\varrho^{\ast}_\theta(\beta)=\cases{ %
(\sqrt{1-\theta}-\sqrt{1-\beta-\theta/2})^2, &\quad $\mbox{if $1/2<
\beta \leq (3-\theta)/4$};$
\vspace*{2pt}\cr
\beta-1/2, & \quad$\mbox{if $(3-\theta)/4<\beta\leq3/4$};$
\vspace*{2pt}\cr
(1-\sqrt{1-\beta})^2, & \quad$\mbox{if $3/4<\beta<1$},$}
\end{eqnarray*}
which is clearly inferior to $\varrho^{\ast}(\beta)$. The boundary
$\varrho^{\ast}_{\theta}(\beta)$ is the same as the one in \citet{Delaigle} based on the marginal $t$-statistics, whereas our
result is based on the $z$-statistics. The $t$-statistic formulation
reduces the demand on the tails of the distributions as shown in
\citet{Delaigle}.
We note that if $\theta\ge1$, Theorem \ref{detect-upper-bound} remains so that the Gaussian
detection boundary is still valid.

%Despite the possible non-Gaussian data, the detection boundary $r=
%is due to the asymptotic normality of each marginal means $\bar{X}_j$.
%Another is that the thresholding statistics are asymptotically
%uncorrelated as %implied from Proposition 1.
% However, we are not saying that the detection boundary is optimal for
%the family of distributions with possible column-wise dependence that
%we are considering in (C2) %and (C3). \footnote{This paragraph is a
%little repeat from the discussion in Section 6, should we delete this
%part?}

%s4 #&#
\section{Power comparison}\label{sec4}

We compare the power of the maximal $L_2$-thresh\-olding test with those
of the HC test and the maximal $L_1$-thresholding test in this section.
Let us first introduce these two tests.

The HC test is based on
%
%e4.1 #&#
\begin{equation}
\label{HCtest} \hat{\mathcal{T}}_{0,n}(s)=\frac{{T}_{0 n}(s)-2p\bar{\Phi}(\lambda
^{1/2}_p(s))}{\sqrt{2p\bar{\Phi}(\lambda^{1/2}_p(s))(1-2\bar{\Phi
}(\lambda^{1/2}_p(s)))}},
\end{equation}
where $T_{0 n} (s)=\sum_{j=1}^p I(Y_{j,n}\geq\lambda_p(s))$. Like
\citet{Delaigle2009}, we consider here a two-sided HC test instead
of a one-sided test treated in \citet{DonohoJin}. With the same
reasoning as Donoho and Jin [(\citeyear{DonohoJin}), page~968], we define the HC test statistic
\[
\hat{\mathcal{M}}_{0 n} = \max_{s \in{\mathcal{S}}} \hat {\mathcal
{T}}_{0,n} (s),
\]
where $\mathcal{S} = (0, 1-\eta]$ for an arbitrary small $\eta$ and is
the same as the maximal $L_2$-thresholding statistic.
Using the same argument for the maximal $L_2$-thresholding statistic,
it can be shown that $\hat{\mathcal{M}}_{0 n}$ attains its maximum
value on $\mathcal{S}_{n}$ given in (\ref{eq:Sn}) as well.

According to \citet{DonohoJin}, under $H_0$,
\[
P\bigl(a(\log p)\hat{\mathcal{M}}_{0 n} -b(\log p,\eta)\leq x\bigr)\to
\exp\bigl(-e^{-x}\bigr),
\]
with the same normalizing sequences as those in Theorem \ref{asy-gumbel}. %$\{a_n(\log
%p)\}$ and $\{b_n(\log p, \eta)\}$.
Let $\mathcal{B}_{\alpha}$ be the same as that of the maximal
$L_2$-thresholding test given in (\ref{eq:L2test}). An $\alpha$ level
HC test rejects~$H_0$ if
%
%e4.2 #&#
\begin{equation}
\hat{\mathcal{M}}_{0 n} >\mathcal{B}_{\alpha}.\label{eq:HCtest}
\end{equation}
%
%=\mathcal{E}_\alpha/a_n(p_s)+b_n(p_s)$ and $\mathcal{B}_{HC,\alpha}$
%is the same as $\mathcal{B}_{\alpha}$. Notice that the set $
%Jin (2004), however, for the fairness of comparison, we keep the same
%set as the maximal $L_2$ test.

Let us introduce the maximal $L_1$-thresholding test statistic. Recall that
\[
T_{1n}(s)=\sum_{j=1}^p|
\sqrt{n}\bar{X}_{j}|I\bigl(|\bar{X}_{j}|>\sqrt {\lambda
_p(s)/n}\bigr).
\]
It can be shown that the mean and variance of {\bf$T_{1n}(s)$} under
$H_0$ are, respectively,
\begin{eqnarray*}
\mu_{T_{1n},0}(s)&=&\sqrt{2/\pi}p^{1-s}\bigl\{1+o(1)\bigr\}
\quad\mbox{and}\\
 \sigma^2_{T_{1n},0}(s)&=&\bigl\{2p^{1-s}
\sqrt{(s/\pi)\log p}\bigr\} \bigl\{1+o(1)\bigr\}.
\end{eqnarray*}

Define
\[
\hat{\mathcal{T}}_{1,n}(s)=\frac{T_{1n}(s)-\hat{\mu
}_{T_{1n},0}(s)}{\tilde{\sigma}_{T_{1n},0}(s)},
\]
where $\hat{\mu}_{T_{1n},0}(s)$ is a sufficiently accurate estimator of
${\mu}_{T_{1n},0}(s)$ in a similar sense to~(\ref{eq:cri1}) and
$\tilde
{\sigma}_{T_{1n},0}^2(s) =2 p^{1-s} \sqrt{(s/\pi)\log p}$. The maximal
$L_1$-thresh\-olding statistic is
\[
\hat{\mathcal{M}}_{1n}=\max_{s\in\mathcal{S}}\hat{
\mathcal{T}}_{1,n}(s),
\]
where, again, $\mathcal{S} = (0, 1-\eta]$.
It can be shown that $\hat{\mathcal{M}}_{1n}=\max_{s\in\mathcal{S}_n}
\hat{\mathcal{T}}_{1,n}(s)$ for the same $S_n$ in (\ref{eq:Sn}).
% where $\mathcal{S}_n$ is defined as the same set as that in maximal
%$L_2$ test.

Using a similar approach to that in Theorem \ref{asy-gumbel}, we can show that
\[
P\bigl(a(\log p)\hat{\mathcal{M}}_{1n}-b(\log p,\eta)\leq x\bigr)\to
\exp\bigl(-e^{-x}\bigr).
\]
Hence, an $\alpha$-level maximal $L_1$-thresholding test rejects the
$H_0$ if
%
%e4.3 #&#
\begin{equation}
\hat{\mathcal{M}}_{1n}>\mathcal{B}_{\alpha}. % =\mathcal{E}_
\label{eq:L1test}
\end{equation}

From (\ref{eq:L2test}), (\ref{eq:HCtest}) and (\ref{eq:L1test}), the
three tests have the same critical values $\mathcal{B}_\alpha$ at
nominal level $\alpha$. This brings convenience for the power
comparison. % among the three tests.
Let us define the power of the three tests %as functions of $(r, \beta)$
\[
%Note that the power functions of the maximal $L_1, L_2$ test are,
%respectively, $
\Omega_{\gamma}(r,\beta):=P(\hat{
\mathcal{M}}_{\gamma n}>\mathcal {B}_{\alpha})
\]
for $\gamma= 0, 1$ and $2$, respectively.
% \mbox{and} \Omega_2(r,\beta):=P(\mathcal{M}_{2n}>\mathcal{B}_{
%and the power function of the HC test is $\Omega_{0}(r,\beta):=P({
%Note that $\Omega_{\gamma}(r,\beta)=P(\mathcal{M}_{\gamma n}>
%for $\gamma=0, 1$ and $2$ where $\mathcal{M}_{\gamma n}=\max_{s\in
Notice that
%
%e4.4 #&#
\begin{equation}
\label{th4decom} \hat{\mathcal{M}}_{\gamma n}=\max_{s\in\mathcal{S}_n}
\bigl\{ \mathcal {T}_{\gamma n}(s)\tilde{e}_\gamma(s)+{\tilde{\sigma
}^{-1}_{T_{\gamma
n}, 0}(s)} \bigl(\mu_{T_{\gamma n},0}(s)-\hat{
\mu}_{T_{\gamma
n},0}(s) \bigr) \bigr\},
\end{equation}
where $\tilde{e}_\gamma(s)={\sigma_{T_{\gamma n}, 0}(s)}/{\tilde
{\sigma
}_{T_{\gamma n}, 0}(s)}$ and
\[
\mathcal{T}_{\gamma n}(s)={{\sigma}^{-1}_{T_{\gamma n}, 0}(s)} {
\bigl(T_{\gamma n}(s)-\mu_{T_{\gamma n},0}(s) \bigr)}=\mathcal{T}_{\gamma
n,1}(s)R_\gamma(s)+
\Delta_{\gamma,0}(s;r,\beta),
\]
in which $R_\gamma(s)={\sigma_{T_{\gamma n}, 1}(s)}/{{\sigma
}_{T_{\gamma n}, 0}(s)}$,
$\mathcal{T}_{\gamma n,1}(s)={{\sigma}^{-1}_{T_{\gamma n}, 1}(s)}
(T_{\gamma n}(s)-\mu_{T_{\gamma n},1}(s) )$ and $\Delta_{\gamma
,0}(s;r,\beta)={{\sigma}^{-1}_{T_{\gamma n}, 0}(s)} (\mu
_{T_{\gamma
n},1}(s)-\mu_{T_{\gamma n},0}(s) )$.
As shown in (\ref{L2-100}), (\ref{HC-100}) and (\ref{L1-100}) in the
\hyperref[app]{Appendix},
\begin{eqnarray*}
 \Delta_{0,0}(s;r,\beta)&=&(s\pi\log p)^{{1}/{4}}p^{1/2-\beta
+s/2}I(r>s)\\
&&{}+L_p^{(6)}p^{1/2-\beta-(\sqrt{s}-\sqrt{r})^2+s/2}I(r<s),
\\
\Delta_{1,0}(s;r,\beta)&=&(s\pi\log p)^{{1}/{4}}(r/s)^{{1}/{4}}
p^{1/2-\beta+s/2}I(r>s)\\
&&{}+L_p^{(6)}p^{1/2-\beta-(\sqrt{s}-\sqrt
{r})^2+s/2}I(r<s)
\end{eqnarray*}
and
\begin{eqnarray*}
 \Delta_{2,0}(s;r,\beta)&=&(s\pi\log p)^{{1}/{4}}(r/s)p^{1/2-\beta
+s/2}I(r>s)\\
&&{}+L_p^{(6)}p^{1/2-\beta-(\sqrt{s}-\sqrt{r})^2+s/2}I(r<
s),
\end{eqnarray*}
where $L_p^{(6)}=\{2(\sqrt{s}-\sqrt{r})\}^{-1}s^{1/4}(\pi\log p)^{-1/4}$.

{Derivations given in the proof of Theorem \ref{th4} in the \hyperref[app]{Appendix} show that}
for $\gamma=0,1$ and $2$,
%
%e4.5 #&#
\begin{equation}
\hat{\mathcal{M}}_{\gamma n}\sim\max_{s\in\mathcal{S}_n}\Delta
_{\gamma
,0}(s;r,\beta), \label{eq:equvi}
\end{equation}
where ``$a\sim b$'' means that the $a/b=1+o_p(1)$. This implies that we
only need to compare $\max_{s\in\mathcal{S}_n}\Delta_{\gamma
,0}(s;r,\beta)$ in the power comparison.

From the established expressions of $\Delta_{\gamma, 0}(s; r,\beta)$,
we note two facts. One is that
if $r>2\beta-1$, for any $s\in(2\beta-1,r)$,
%
%e4.6 #&#
\begin{eqnarray} \label{eq:4.20a}
\Delta_{2,0}(s;r,\beta)/\Delta_{1,0}(s;r,
\beta)&=&(r/s)^{3/4}>1\quad \mbox{and}
\nonumber
\\[-8pt]
\\[-8pt]
\nonumber
\Delta_{1,0}(s;r,\beta)/
\Delta_{0,0}(s;r,\beta)&=&(r/s)^{1/4}>1.
\end{eqnarray}
%
%and the three $\Delta_{\gamma,0}$s are asymptotically equivalent for
%any $s>r$. %$\Delta_{2,0},\Delta_{1,0}$ and $\Delta_{HC,0}$ diverge to
%infinity %as $n\to\infty$.
The other is if $r\in(\varrho^*(\beta),2\beta-1]$, asymptotically,
%
%e4.7 #&#
\begin{equation}
\Delta_{0,0}(s;r,\beta)=\Delta_{1,0}(s;r,\beta)=\Delta
_{2,0}(s;r,\beta)\qquad \mbox{for all $s\in\mathcal{S}$}. \label{eq:4.20b}
\end{equation}
Hence, when $(r,\beta)$ lies just above the detection boundary, the
three $\Delta_{\gamma, 0}$ functions are the same.
If $(r, \beta)$ moves further away from the detection boundary so that
$r > 2 \beta-1$, there will be a clear ordering among the $\Delta
_{\gamma,0}$ functions.
%The difference among the $\Delta_{\gamma,0}$ will be more substantial
%if $r > 1-\eta$, where $1-\eta$ is the maximum value of $s$.
The following theorem summarizes the relative power performance.

%th4 #&#
\begin{theorem}
\label{th4} Assume \textup{(C.1)--(C.5)} and (\ref{eq:cri1}) hold. For any
given significant level $\alpha\in(0,1)$, the powers of the HC, the
maximal $L_1$ and $L_2$-thresholding tests under $H_1$ as specified in
\textup{(C.4)} satisfy, as $n \to\infty$,
%fixed quantities, it seems strange when we saying with probability one
%this happen. On the other hand, %indeed, when we compare the powers,
%we know with probability one, $\hat{\mathcal{M}}_{0 n}<\hat{
%we %should have some randomness here. What should we say here?}} }
%e4.8 #&#
\begin{equation}
\label{pow-order-1} \Omega_{0}(r,\beta)\leq\Omega_1(r,\beta)
\leq\Omega_2(r,\beta)\qquad \mbox{for $r>2\beta-1$}%\mbox{when $r>\varrho^*(\beta)$}.
\end{equation}
and $\Omega_{\gamma}(r,\beta)$s are asymptotic equivalent for $r\in
(\varrho^*(\beta),2\beta-1]$.
\end{theorem}

{The theorem indicates that when $(r,\beta)$ is well above the
detection boundary such that $r > 2 \beta-1$, there is a clear
ordering in the power among the three tests, with the $L_2$ being the
most powerful followed by the $L_1$ test. However, when $(r,\beta)$ is
just above the detection boundary such that $r\in(\varrho^*(\beta
),2\beta-1]$, the three tests have asymptotically equivalent powers. In
the latter case, comparing the second order terms of $\hat{\mathcal
{M}}_{\gamma n}$
may lead to differentiations among the powers of the three tests.
However, it is a rather technical undertaking to assess the impacts of
the second order terms.} %However, the
%
%hold for $\max_{s\in\mathcal{S}_n}\Delta_{\gamma, 0}(s; r, \beta)$
%with probability one when $r > 2 \beta-1$. Therefore, we do not need
%to evaluate the small order term %related to $\mathcal{T}_{\gamma
%n,1}(s){R}_\gamma(s)$. However, if $r\in(\varrho^*(\beta),2\beta-1]$, $
%%asymptotic equivalent. In such case, the second order terms will
%affect the power comparison, which is difficult to evaluate. Since $
%dominates when $r\in(\varrho^*(\beta),2\beta-1]$, the power functions $
%to 1 at the same rate.}
% the corresponding strict ordering on $\max_{s\in\mathcal{S}}\Delta_{
%samples, the ordering in (\ref{eq:4.20a}) may also produce an ordering
%in the power as conveyed in our simulation study reported in the next
%section.
The analysis conducted in Theorem \ref{th4} is applicable to the setting of
Gaussian data with $n=1$ and $\bolds{\Sigma}$ satisfying (C.3), which is
the setting commonly assumed in the investigation of the detection
boundary for the HC test [\citet{DonohoJin}; \citet{HallJin2010}
and Arias-Castro, Bubeck and
Lugosi (\citeyear{Arias-Castro2012a})].
Specifically, the power ordering among the three maximal thresholding
tests in Theorem \ref{th4} remains but under lesser conditions (C.3)--(C.5).
Condition~(C.1) is not needed since the Gaussian assumption allows us to
translate the problem to $n=1$ since the sample mean is sufficient.
Condition (C.2) is automatically satisfied for the Gaussian distribution.
The condition (\ref{eq:cri1}) is met for the Gaussian data, as we have
discussed in Section~\ref{sec2}.

%s5 #&#
\section{Simulation results}\label{sec5}

We report results from simulation experiments which were designed to
evaluate the performance of the maximal $L_1$ and $L_2$-thresholding
tests and the HC test.
% as well as the test based on the false discovery rate (Benjamini and
%Hochberg, 1995).
The purpose of the simulation study is to confirm the
theoretical findings that there is an
ordering in the power among the three tests discovered in Theorem \ref{th4}.
%that the maximal $L_2$ test is better that the maximal $L_1$, which
%is better than the HC test.

Independent and identically distributed $p$-dim random vectors $\mathbf
{X}_i$ were
generated according to
\[
\mathbf{X}_i=\mathbf{W}_i+\bolds{\mu},\qquad i=1,\ldots,n,
\]
where $\mathbf{W}_i=(W_{i1},\ldots,W_{ip})^T$ is a stationary random
vector and $\{W_{ij}\}_{j=1}^p$ have the same marginal distribution $F$.
In the simulation, $\mathbf{W}_i$ was generated from a $p$-dimensional
multivariate Gaussian distribution with zero mean and covariance
$\bolds{\Sigma} = (\sigma_{ij})_{p \times p}$, where $\sigma_{ij} =
\rho^{|i-j|}$ for $\rho=0.3$ and $0.5$, respectively.
%To introduce dependence among $\{W_{ij}\}_{j=1}^p$, we employed the
%$p$-dimensional Gaussian copula $C_{Gauss, 0, V}$, that corresponds to
%$p$-dimensional normal distribution
%The distribution we experimented for $\mathbf{W}_i$ %= (W_{i1},
%was $$C_{Gauss, 0, V}\{\Phi(w_1), \cdots, \Phi(w_p)\}$$
%where $\Phi$ is the standard normal distribution. This is a bona fide
%multivariate Gaussian distribution.
%The second distribution for $\mathbf{W}_i$ was
%$$C_{Gauss, 0, V}\{F_{t_3}(w_1), \cdots, F_{t_3}(w_p)\}$$
%where $F_{t_3}$ is the standardized $t_3$ distribution such that it
%has unit variance. Hence, both distribution had the same dependence
%structure as dictated by the Gaussian copula, but different marginal
%distributions.

The simulation design on $\bolds{\mu}$ had the sparsity parameter
$\beta=0.6, 0.7$ and $0.8$, respectively, and the signal strength $r
=0.1, 0.3, 0.5, 0.6, 0.8, 0.9, 1.1$ and $1.2$, respectively. We chose
two scenarios on the dimension and sample size combinations: (a) a
large $p$, small $n$ setting
and (b) both $p$ and $n$ are moderately large. For scenario (a), we
chose $p=\exp(c_0n^{0.3}+c_1)$, where $c_0=1.90$ and $c_1=2.30$ so that
the dimensions $p$ were 2000 and 20,000, and the sample sizes $n$ were
$30$ and 100, respectively. We note that under the setting $\beta=0.8$,
there were only $4$ and 7 nonzero means, respectively, among the 2000
and 20,000 dimensions. And those for $\beta=0.7$ were $9$ and $19$,
respectively, and those for $\beta=0.6$ were $20$ and $52$,
respectively. These were quite sparse. For scenario (b), we chose
$p=n^{1.25}+184$ such that $(p,n)=(500,100)$ and $(p,n)=(936, 200)$.

The maximal $L_2$-test statistic $\hat{\mathcal{M}}_{2n}$ was
constructed using $\tilde{\mu}_{T_{2n},0}(s)$ and $\tilde{\sigma
}_{T_{2n}, 0}(s)$ given in (\ref{eq:meanTn0}) and (\ref{eq:varTn0}),
respectively, as the mean and standard deviation estimators.
The maximal $L_1$ test statistic and the HC test statistic, $\hat
{\mathcal{M}}_{1n}$ and $\hat{\mathcal{M}}_{0 n} $, were constructed
similarly using the leading order mean and standard deviation under
$H_0$. The set of thresholding level $\mathcal{S}$ was chosen to be
$(0, 1-\eta]$ with $\eta=0.05$.
%Since the empirical relative performance of the three tests with
%different marginal distributions were largely the same, we only report
%the simulation results for the first model with the Gaussian marginals.

Figures~\ref{figure1}--\ref{figure4} display the average empirical sizes
and powers of the HC, the maximal $L_1$ and $L_2$-thresholding tests
based on 20,000 simulations, with
Figures~\ref{figure1}--\ref{figure2} for scenario (a) and Figures~\ref{figure3}--\ref{figure4} for scenario (b).} % under Gaussian
%distribution.}
%Tables 1 and 2 for the Gaussian data and Table 3 and 4 for the
%standardized $t_3$.
%and the FDR tests, where table 1 and 3 contains the result for Gassian
%marginal case and table 2 and 4 are for $t_3$ distribution marginal
%case.
%Each entry in the tables was based on 20,000 replications.
To make the power comparison fair and conclusive, we adjusted the
nominal level of the tests so that the simulated sizes of the tests
were all around
$\alpha=0.05$, with the HC having slightly larger sizes than those of
the maximal $L_1$ test, and the sizes of the maximal $L_1$ test were
slightly larger than those of the maximal $L_2$ test. These were
designed to rule out potential ``favoritism'' in the power comparison
due to advantages in the sizes of the maximal $L_2$ and/or $L_1$ tests.
%%Because the empirical sizes appeared differently, adjusted sizes $
%and close to $\alpha=0.05$. The powers of the tests were compared at
%the adjusted sizes $\alpha^*$.

%f1 #&#
\begin{figure}

\includegraphics{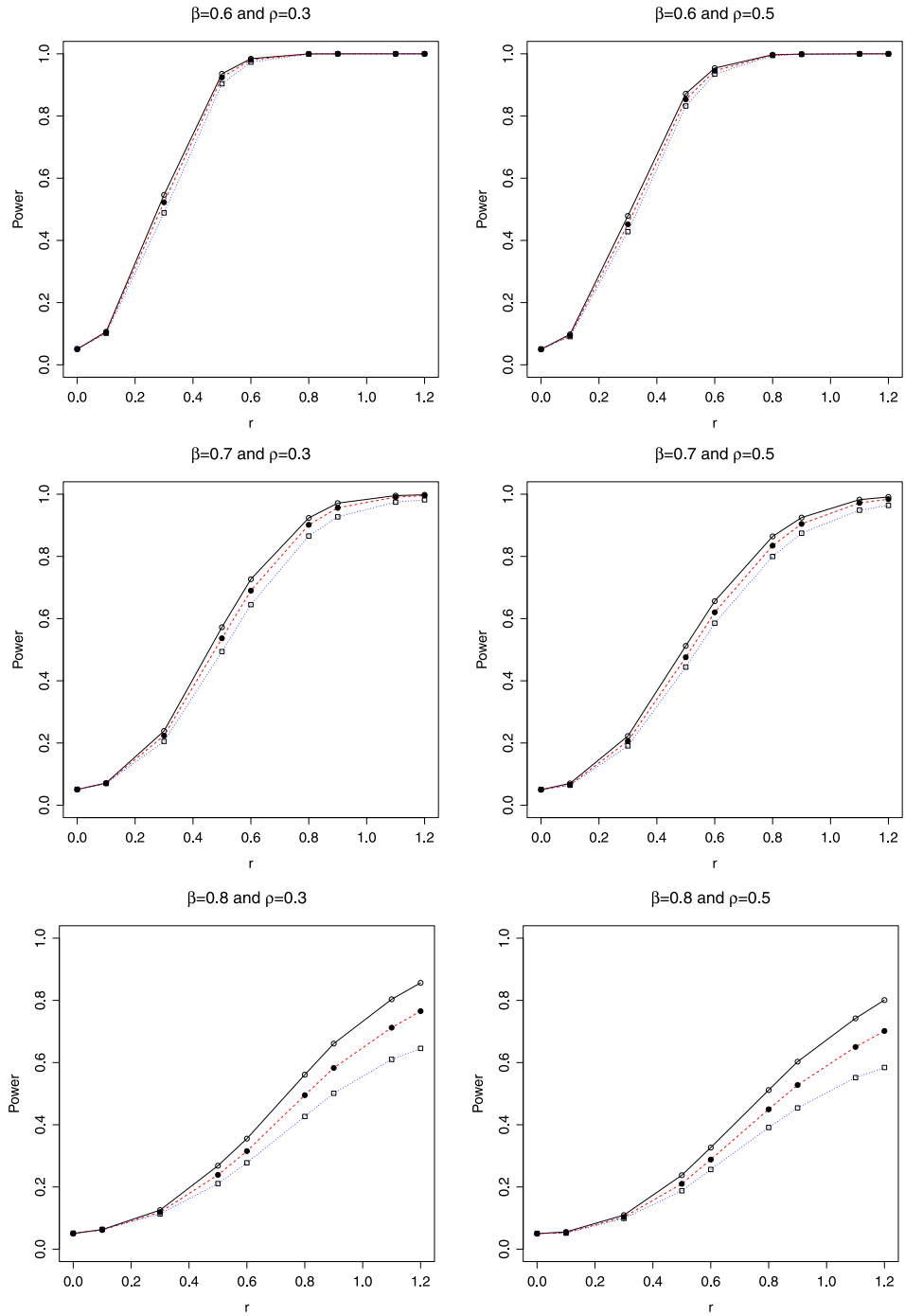}

\caption{Empirical sizes and powers of the HC (dotted lines with
squares), the maximal $L_1$- (dashed lines with dots) and $L_2$- (solid
lines with circles) thresholding tests when $p=2000$ and $n=30$ with
the marginal distribution the standard normal.}\label{figure1}
\end{figure}

%f2 #&#
\begin{figure}

\includegraphics{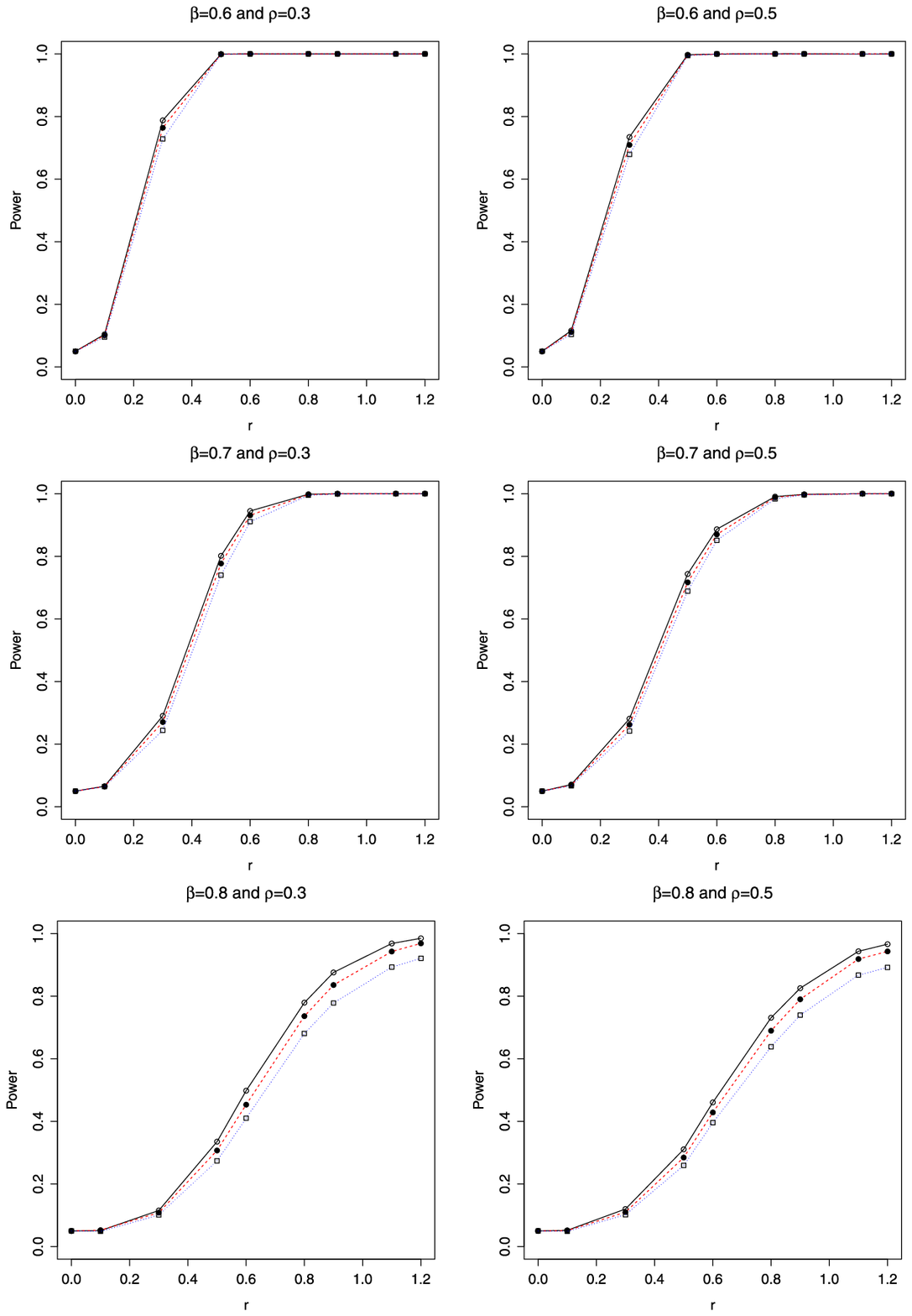}

\caption{Empirical sizes and powers of the HC (dotted lines with
squares), the maximal $L_1$- (dashed lines with dots) and $L_2$- (solid
lines with circles) thresholding tests when $p=20\mbox{,}000$ and $n=100$ with
the marginal distribution the standard normal.}\label{figure2}
\end{figure}

%f3 #&#
\begin{figure}

\includegraphics{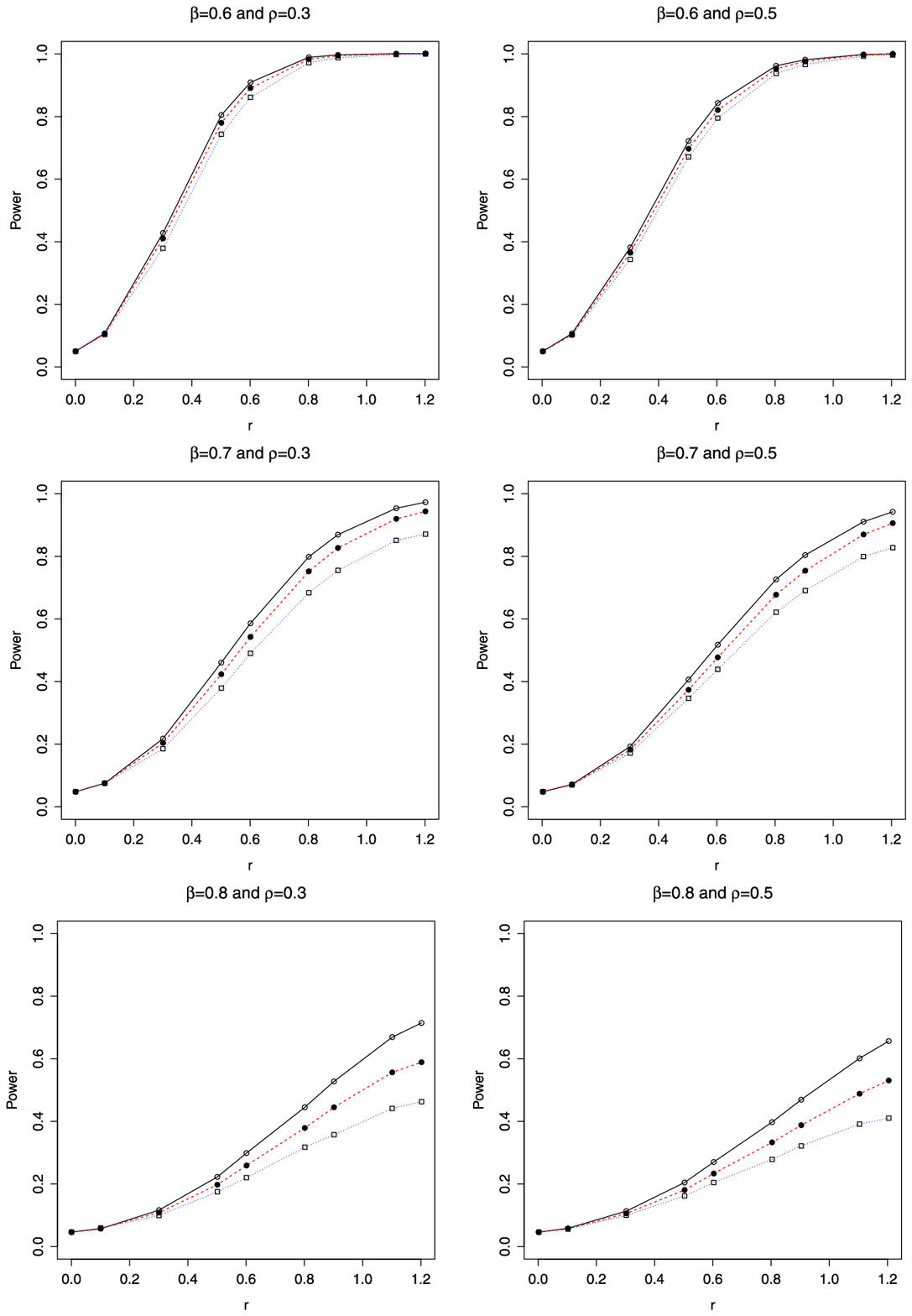}

\caption{Empirical sizes and powers of the HC (dotted lines with
squares), the maximal $L_1$- (dashed lines with dots) and $L_2$- (solid
lines with circles) thresholding tests when $p=500$ and $n=100$ with
the marginal distribution the standard normal.}\label{figure3}
\end{figure}

%f4 #&#
\begin{figure}

\includegraphics{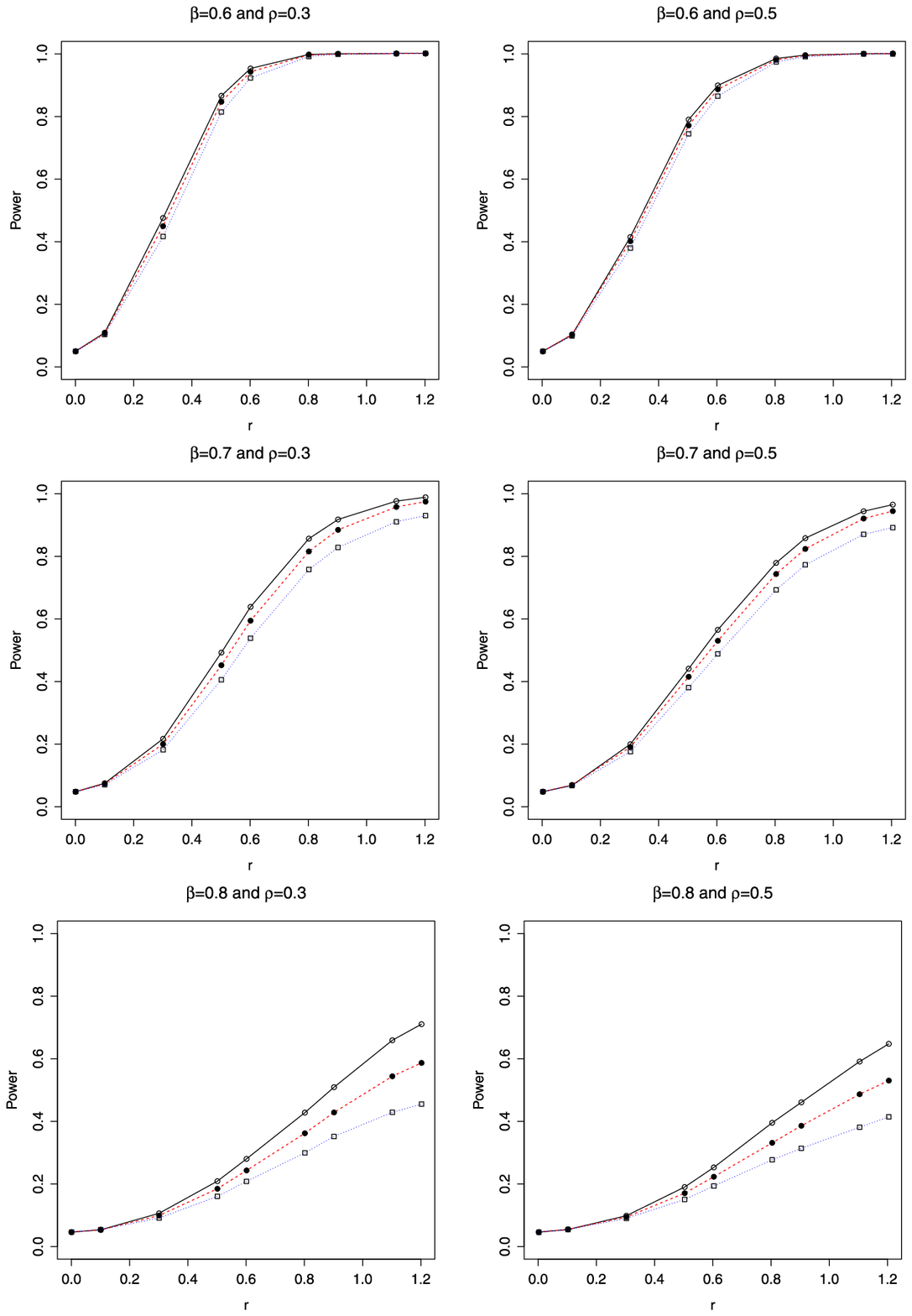}

\caption{Empirical sizes and powers of the HC (dotted lines with
squares), the maximal $L_1$- (dashed lines with dots) and $L_2$- (solid
lines with circles) thresholding tests when $p=936$ and $n=200$ with
the marginal distribution the standard normal.}\label{figure4}
\end{figure}

Figures~\ref{figure1}--\ref{figure4} show that the power of the tests
were the most influenced by the signal strength parameter $r$, followed
by the sparsity $\beta$. The powers were insensitive to the level of
dependence $\rho$, which confirmed our finding that the thresholding
largely removes the dependence. The observed ordering in the empirical
power shown in Figures~\ref{figure1}--\ref{figure4} were consistent to
the conclusions in Theorem \ref{th4}. We observed that in all the simulation
settings, despite some size advantages by the HC test and/or the
maximal $L_1$ test, the maximal $L_2$ test had better power than the
maximal $L_1$ and the HC test, and the maximal $L_1$ test had better
power than the HC test. We find that for each fixed level of sparsity
$\beta$, when the signal strength $r$ was increased so that $(r, \beta
)$ moved away from the detection boundary $r=\varrho^{\ast}(\beta)$,
the difference among the powers of the three tests was enlarged. This
was especially the case for the most sparse case of $\beta=0.8$ and was
indeed confirmatory to Theorem \ref{th4}.
The simulated powers of the three tests were very much the same at
$r=0.1$ and were barely changed even when both $n$ and $p$ were
increased. This was consistent with the fact that $r=0.1$ is below the
detection boundary for $\beta=0.7$ and 0.8 considered in the simulation.

%s6 #&#
\section{Discussion}\label{sec6}

Our analysis shows that there are alternative $L_1$ and $L_2$
formulations to the HC test which attain the detection boundary
$r=\varrho^{\ast}(\beta)$ of the HC test. The tests based on the $L_1$
and $L_2$ formulations are more powerful than the HC test when the $(r,
\beta)$ pair is away from the detection boundary such that $r > 2\beta-1$.
The three tests have asymptotically equivalent power when $(r,\beta)$
is just above the detection boundary.

The detection boundary $r=\varrho^{\ast}(\beta)$ coincides with that of
the HC test discovered in \citet{DonohoJin} for the Gaussian data
with independent components. That the three tests considered in this
paper attain the detection boundary $r=\varrho^{\ast}(\beta)$ under the
considered sub-Gaussian setting with column-wise dependence can be
understood in two aspects. One is that the three test statistics are
all directly formulated via the marginal sample means $\bar{X}_j$ which
are asymptotically normally distributed; the other is that the
thresholding statistics are asymptotically uncorrelated as implied from
Proposition \ref{chap4-cor2}.

According to \citet{Ingster} and \citet{DonohoJin}, $r=\varrho
^{\ast
}(\beta)$ is the optimal detection boundary for Gaussian distributed
data with independent components. However, it may not be optimal for
the dependent nonparametric setting considered in this paper. Indeed,
for weakly dependent Gaussian data, \citet{HallJin2010} showed that the
detection boundary $r=\varrho^{\ast}(\beta)$ can be lowered by
utilizing the dependence. The latter was carried out by
pre-transforming the data with~$\mathbf{L}$, the inverse of the
Cholesky decomposition of $\bolds{\Sigma}$, or an empirical estimate
of $\mathbf{L}$ and then conducting the HC test based on the
transformed data. It is expected that the main results of this paper on
the relative performance of the three tests would remain valid for the
transformed data. \citet{HallJin} and \citet{Delaigle2009}
studied the detection boundary for dependent data and \citet{Cai2012}
studied the boundary for detecting mixtures with a general known
distribution. However, the optimal detection boundary under the
dependent sub-Gaussian distribution setting is still an open problem.

%sA #&#
\begin{appendix}
%sA #&#
\section*{Appendix: Technical details}\label{app}

In this Appendix we provide proofs to Theorems \ref{asy-gumbel}, \ref{detect-upper-bound}
and \ref{th4} reported in
Sections~\ref{sec3} and \ref{sec4}.
Throughout this Appendix we use $L_p=C\log^b(p)$ to denote slow varying
functions
for some constant $b$ and positive constant $C$, and $\phi(\cdot)$ and
$\bar{\Phi}(\cdot)$ for the density and survival functions of the
standard normal distribution, respectively. Let $\rho_k$ be the
correlation coefficient between $W_{i1}$ and $W_{i(k+1)}$, and
write $\rho_1=\rho$ for simplicity and $\mu_j=E(X_{ij})$ for $i\in\{
1,\ldots, n\}$ and $j\in\{1,\ldots,p\}$. Put $\lambda_p(s)=2 s \log p$.

\begin{pf*}{Proof of Theorem \ref{asy-gumbel}}
%By Theorem 1 in Galambos (1973), we only need to check the asymptotic
%independence of $\{\mathcal{T}_2(s_i): s_i\in\mathcal{S}_n\}$ under
%the null hypothesis.
%This could be checked if we can show that $\mathcal{T}_{2}(s)$
%converges uniformly on $\mathcal{S}$ to a zero mean Gaussian process $
%show this, we divide the proof into two parts. Part (i) established
%the asymptotic normality of $\mathcal{T}_{2}(s)$ at finite points $
%the $\mathcal{T}_{2}(s)$.
Let $u=\bar{\Phi}(\lambda^{1/2}_p(s))$. Write $\mathcal
{J}_2(u):=\hat
{\mathcal{T}}_{2,n}(s)$ and
\[
\mathcal{M}_{2n}=\max_{s\in(0,1-\eta]}\hat{\mathcal
{T}}_{2,n}(s)=\max_{u\in[u_0, 1/2)}\mathcal{J}_2(u),
\]
where $u_0=\bar{\Phi}(\lambda_p^{1/2}(1-\eta))$. Using the same
technique for the proof of Theorem 1 in Zhong, Chen and Xu (\citeyear{ZhenChenXu}), it may be
shown that the joint asymptotic normality of $\mathcal{T}_{2,n}(s)$ at
any finite points $\underline{s}=(s_1,\ldots,s_d)^T$. This is
equivalent to the joint asymptotic normality of $\mathcal{J}_{2}(u)$ at
$u_i=\bar{\Phi}(\sqrt{2s_i\log p})$ for $i=1, \ldots, d$.

We want to show the tightness of the process $\mathcal{J}_2(u)$. Let
$f_{n,u}(x)=\sigma_0^{-1}(u)\times x^2I\{|x|>g(u)\}$, where $g(u)=\bar{\Phi
}^{-1}(u)$, $\sigma_0^2(u)=\sigma_0^2(p;s)$ and $\sigma
_0^2(p;s)=\break \sigma
_{T_{2n},0}^2(s)/p$. Write
\[
\mathcal{J}_2(u)=p^{-1/2}\sum_{j=1}^p
\bigl\{f_{n,u}\bigl(|\sqrt{n}\bar {X}_j|\bigr)-E
\bigl(f_{n,u}\bigl(|\sqrt{n}\bar{X}_j|\bigr)\bigr) \bigr\}.
\]
%
%$p\sigma^2_0(p;s)=\sigma_{T_{2n},0}^2(s)$ defined in (\ref{eq:varTn0}).
Based on the finite dimensional convergence of $\mathcal{J}_2(u)$ and
Theorem 1.5.6 in Van der Vaart and Wellner (\citeyear{Wellner}), we only need to show
the asymptotically equicontinuous of $\mathcal{J}_2(u)$, that is, for
any $\varepsilon>0$ and $\eta>0$ there exists a finite partition
$\Lambda=\bigcup_{i=1}^k \Lambda_i$ such that
%
%eA.1 #&#
\begin{equation}
\label{equaconti} \lim\sup_{n\to\infty}P^*\Bigl\{\max
_{1\leq i\leq k}\sup_{u,v\in
\Lambda
_i}\bigl|\mathcal{J}_2(u)-
\mathcal{J}_2(v)\bigr|>\varepsilon\Bigr\}<\eta,
\end{equation}
where $P^*$ is the outer probability measure.

Define $
\mathscr{F}_n=\{f_{n,u}(|\sqrt{n}\bar{X}_j|)=\sigma_0^{-1}(u)|\sqrt {n}\bar{X}_j|^2I\{|\sqrt{n}\bar{X}_j|>g(u)\}\dvtx u\in\Lambda
:=[u_0,1/2)\}
$
and $\rho(f_{n,u}-f_{n,v})=[E\{f_{n,u}(|\sqrt{n}\bar
{X}_j|)-f_{n,v}(|\sqrt{n}\bar{X}_j|)\}^2]^{1/2}$. It can be shown that
if $u>v$,
\[
\rho(f_{n,u}-f_{n,v})^2=\bigl\{2-2
\sigma_0^{-1}(u)\sigma_0(v)\bigr\} \bigl
\{1+o(1)\bigr\}.
%& +\Phi(\sqrt{2t\log p})-\Phi(\sqrt{2s\log p})\big\}\{1+o(1)\}.
\]
Thus, for every $\delta_n\to0$, $\sup_{|u-v|<\delta_n}\rho
(f_{n,u}-f_{n,v})\to0$, which implies that for each $\delta>0$,
$\Lambda$ can be partitioned into finitely many sets $\Lambda
_1,\ldots
,\Lambda_k$ satisfying
\[
\max_{1\leq i\leq k}\sup_{u,v\in\Lambda_i}\rho
(f_{n,u}-f_{n,v})<\delta.
\]
Let $N_0:=N(\varepsilon,\mathscr{F}_n, \rho)$ be the bracketing number,
the smallest number of functions $f_1,\ldots, f_{N_0}$ in $\mathscr
{F}_n$ such that for each $f$ in $\mathscr{F}_n$ there exists an $f_i$
($i\in\{1,\ldots,{N_0}\}$) satisfying $\rho(f-f_i)\leq\varepsilon
\leq
1$. Applying Theorem 2.2 in \citet{Andrews}, if the
following two conditions hold for an even integer $Q\geq2$ and a real
number $\gamma>0$ such that
%
%eA.2 #&#
%eA.3 #&#
\begin{eqnarray}
\label{condi-1} \sum_{d=1}^\infty
d^{Q-2}\alpha(d)^{{\gamma}/{(Q+\gamma)
}}&<&\infty \quad\mbox{and}
\\
\label{condi-2} \int_{0}^1\varepsilon^{-{\gamma}/{(2+\gamma)}}N(
\varepsilon ,\mathscr {F}_n, \rho)^{1/Q}\,d\varepsilon&<&\infty,
\end{eqnarray}
we have for $n$ large enough
$
\|\sup_{\rho(f_{n,u}-f_{n,v})<\delta\atop u,v\in\Lambda
_i}|\mathcal
{J}_2(u)-\mathcal{J}_2(v)| \|_Q<k^{-1/Q}\eta\varepsilon$.

Invoking the maximal inequality of \citet{Pisier}, it follows that
\[
\Bigl\|\max_{1\leq i\leq
k}\mathop{\sup_{\rho(f_{n,u}-f_{n,v})<\delta}}_{
s,t\in\Lambda_i}\bigr|
\mathcal{J}_2(u)-\mathcal{J}_2(v)\bigr| \Bigr\|_Q<\eta
\varepsilon.
\]
Now using the Markov inequality, we get for $n$ large enough
\begin{eqnarray*}
&&P^*\Bigl\{\max_{1\leq i\leq k}\sup_{u,v\in\Lambda_i}\bigl|\mathcal
{J}_2(u)-\mathcal{J}_2(v)\bigr|>\varepsilon\Bigr\}
\\
&&\qquad \leq \Bigl\|\max_{1\leq i\leq k}\mathop{\sup_{\rho(f_{n,u}-f_{n,v})<\delta
}}_{ u,v\in\Lambda_i}\bigl|
\mathcal{J}_2(u)-\mathcal{J}_2(v)\bigr|\Bigr \| _Q/
\varepsilon< \eta.
\end{eqnarray*}
Hence, the condition (\ref{equaconti}) holds and $\mathcal{J}_2(u)$ is
asymptotically tight.

It remains to show (\ref{condi-1}) and (\ref{condi-2}) hold. For
(\ref
{condi-2}), we note that $\mathscr{F}_n$ is a V-C class for each $n$.
This is because
\[
\mathscr{G}_n=\bigl\{f_{n,u}(x)=\sigma_0^{-1}(u)I
\bigl(x>g(u)\bigr)\dvtx u\in(u_0,1/2)\bigr\}
\]
is a V-C class with VC index 2. Let $\varphi(x)=x^2$. Then $\mathscr
{F}_n=\varphi\cdot\mathscr{G}_n$ is a V-C class by Lemma 2.6.18 in Van
der Vaart and Wellner (\citeyear{Wellner}). Let $G_n(x,u_0)=\sup_{u\in\Lambda
}|f_{n,u}(x)|$ be the envelop function for class $\mathscr{F}_n$.
Clearly, we can take $G_n(x,u_0)=\sigma_0^{-1}(u_0)x^2$. It is easy to
see that $\rho\{G_n(|\sqrt{n}\bar{X}_i|,u_0)\}<\infty$ for a constant
$u_0>0$. Applying a result on covering number of V-C classes [Theorem
2.6.7, Van der Vaart and Wellner (\citeyear{Wellner})], we get $N(\varepsilon
,\mathscr
{F}_n, \rho)\leq K \varepsilon^{-2}$
for a universal constant $K$. It can be verified that if $Q>2+\gamma$,
then (\ref{condi-2}) holds.
The condition~(\ref{condi-1}) follows from the assumption that \mbox{$\rho
_Z(d)\leq C\alpha^d$}.
%$\sum_{d=1}^\infty d^{Q-2}\rho_X(d)^{\frac{\gamma}{Q+\gamma}}<\infty.$

As a result, $\mathcal{J}_2(u)$ converge to a zero mean Gaussian
process $\mathcal{N}_2(u)$ with
%
%eA.4 #&#
\[
\operatorname{Cov}\bigl(\mathcal{N}_2(u),\mathcal{N}_2(v)\bigr)=
\frac{\sigma
_0(u)}{\sigma
_0(v)}=\exp\biggl(-{\frac{1}{2}}\bigl[\log\bigl\{
\sigma_0^2(v)\bigr\}-\log\bigl\{\sigma
_0^2(u)\bigr\} \bigr]\biggr)
\]
for $u<v$. It can be shown that there exists an Ornstein--Uhlenbeck (O--U)
process $\mathcal{U}_2(\cdot)$ with mean zero 0 and $E(\mathcal
{U}_2(u)\mathcal{U}_2(v))=\exp\{-|u-v|\}$ such that $\mathcal
{N}_2(u)=\mathcal{U}_2({\frac{1}{2}}\log\{\sigma_0^2(u)\})$.
Therefore, by a
result for the O--U process in Leadbetter, Lindgren and
Rootz{\'e}n [(\citeyear{Leadbetter}), page 217],
\begin{eqnarray*}
P\Bigl(\max_{s\in\mathcal{S}}\hat{\mathcal{T}}_{2,n}(s)<B_{\tau
_n}(x)
\Bigr)&=&P\Bigl(\max_{u\in\Lambda}\mathcal{N}_2(u)<B_{\tau_n}(x)
\Bigr)\bigl\{ 1+o(1)\bigr\}
\\
&=&P\Bigl(\max_{u\in(0,\tau_n)}\mathcal{U}_2(u)<B_{\tau_n}(x)
\Bigr)\to\exp\bigl\{ -\exp (-x)\bigr\},
\end{eqnarray*}
where $\tau_n={\frac{1}{2}}\log\{{\sigma_0^2({\frac
{1}{2}})}/{\sigma_0^2(u_0)}\}$,
$B_{\tau_n}(x)=(x+b^*(\tau_n))/a(\tau_n)$, $a(t)=\break  (2\log(t))^{1/2}$ and
$b^*(t)=2\log(t)+2^{-1}\log\log(t)-{\frac{1}{2}}\log(\pi)$. From
(\ref
{eq:varTn0}), we have $\tau_n=\frac{1-\eta}{2}\log p\{1+o(1)\}$. Since
\begin{eqnarray*}
a(\tau_n)\max_{u\in(0,\tau_n)}\mathcal{U}_2(u)-b^*(
\tau_n)&=&\frac
{a(\tau
_n)}{a(\log p)}\Bigl[a(\log p)\max_{u\in(0,\tau_n)}
\mathcal {U}_2(u)-b^*(\log p)\Bigr]
\\
&&{} +\frac{a(\tau_n)}{a(\log p)}b^*(\log p)-b^*(\tau_n),
\end{eqnarray*}
${a(\tau_n)}/{a(\log p)}\to1$ and
\begin{eqnarray*}
\frac{a(\tau_n)}{a(\log p)}b^*(\log p)-b^*(\tau_n)&=&\frac{a(\tau
_n)}{a(\log p)}
\bigl[b^*(\log p)-b^*(\tau_n)\bigr]
\\
&&{} +b^*(\tau_n)\biggl[\frac{a(\tau_n)}{a(\log p)}-1\biggr]\to-\log
\frac
{(1-\eta)}{2},
\end{eqnarray*}
we have
\begin{eqnarray*}
&&a(\tau_n)\max_{u\in(0,\tau_n)}\mathcal{U}_2(u)-b^*(
\tau_n)\\
&&\qquad=a(\log p)\max_{u\in(0,\tau_n)}\mathcal{U}_2(u)-
\biggl(b^*(\log p)+\log\frac{(1-\eta)}{2}\biggr).
\end{eqnarray*}
Finally, note that $b^*(\log p)+\log\frac{(1-\eta)}{2}%=2\log(\log
%p)+2^{-1}\log\log(\log p)-2^{-1}\log({4}\pi/{(1-\eta)^2})
=b(\log p, \eta)$. This finishes the proof of Theorem \ref{asy-gumbel}.
\end{pf*}

\begin{pf*}{Proof of Theorem \ref{detect-upper-bound}} %To obtain the detection
%boundary for the maximal $L_2$ test,
(i). The proof is made under four cases. For each case, we find the
corresponding detectable region and the
union of the four regions are the overall detectable region of the
thresholding test. Basically, we show for any $(\beta,r)$ above
$\varrho
^*(\beta)$ within one of the four cases, there exists at least one
threshold level $s$ such that $H_1$ is detectable. For notation
simplification, we only keep the leading order terms for $\mu
_{T_{2n},1}(s)-\mu_{T_{2n},0}(s)$, $\sigma_{T_{2n,1}}(s)$, $\sigma
_{T_{2n,0}}(s)$ and $\Delta_2(s;r,\beta)$.

\textit{Case} 1: $s\leq r$ and $s\leq\beta$. In this case,
$\mu_{T_{2n},1}(s)-\mu_{T_{2n},0}(s)=L_pp^{1-\beta}$ and
$\sigma_{T_{2n},1}(s)=\sigma_{T_{2n},0}(s)=L_pp^{(1-s)/2}$. Hence,
\[
\Delta_2(s;r,\beta)=\frac{\mu_{T_{2n},1}(s)-\mu
_{T_{2n},0}(s)}{\sigma
_{T_{2n},1}(s)}=L_pp^{(1+s-2\beta)/2}.
\]
So to make %the test detectable, i.e. such that
$(\mu_{T_{2n},1}(s)-\mu_{T_{2n},0}(s))/\sigma_{T_{2n},1}(s)\to
\infty$,
we need $s> 2\beta-1$.
It follows that the detectable region for
this case is $r\geq2\beta-1$. Specifically, if we
select $s=\min\{r,\beta\}$, we arrive at the best divergence rate for
$\Delta_2(s;r,\beta)$ of order $L_pp^{(1+\min\{r,\beta\}-2\beta)/2}$.

\textit{Case} 2: $s\leq r$ and $s>\beta$. In this case,
$\mu_{T_{2n},1}(s)-\mu_{T_{2n},0}(s)=L_pp^{1-\beta}$,
$\sigma_{T_{2n},1}(s)=L_pp^{(1-\beta)/2}$, and
$\sigma_{T_{2n},0}(s)=L_pp^{(1-s)/2}$. Then,
\[
\Delta_2(s;r,\beta)=\frac{\mu_{T_{2n},1}(s)-\mu
_{T_{2n},0}(s)}{\sigma
_{T_{2n},1}(s)}=L_pp^{(1-\beta)/2}.
\]
So the detectable region in the $(\beta,r)$ plane is $r>\beta$. In
this region, the best divergence rate of $\Delta_2$ is of order
$L_pp^{(1-\beta)/2}$ for any $\beta<s\leq r$.

\textit{Case} 3: $s> r$ and $s\leq(\sqrt{s}-\sqrt{r})^2+\beta$.
The case is equivalent to $\sqrt{r}< \sqrt{s}\leq
(r+\beta)/(2\sqrt{r})$ and
$\mu_{T_{2n},1}(s)-\mu_{T_{2n},0}(s)=L_pp^{1-(\sqrt{s}-\sqrt
{r})^2-\beta}$,
$\sigma_{T_{2n},1}(s)=\sigma_{T_{2n},0}=L_pp^{(1-s)/2}$. Then
%
%eA.5 #&#
\begin{equation}
\label{4:eq:3} \Delta_2(s;r,\beta)=\frac{\mu_{T_{2n},1}(s)-\mu
_{T_{2n},0}(s)}{\sigma
_{T_{2n},1}(s)}=L_pp^{{{1}/{2}}-\beta+r-(\sqrt{s}-2\sqrt{r})^2/2}.
\end{equation}
To ensure (\ref{4:eq:3}) diverging to
infinity, we need
\[
2\sqrt{r}-\sqrt{1-2\beta+2r}<\sqrt{s}<2\sqrt{r}+\sqrt{1-2\beta+2r}.
\]
Thus, the detectable region must satisfy
\begin{eqnarray*}
\sqrt{r}&<&(r+\beta)/(2\sqrt{r}), \qquad 1-2\beta+2r>0\quad \mbox{and}\\
 2\sqrt {r}-\sqrt{1-2
\beta+2r}&\leq&(r+\beta)/(2\sqrt{r}).
\end{eqnarray*}
This translates to
\[
\beta-{\tfrac{1}{2}}<r<\beta\quad\mbox{and } \mbox{either } r\leq\beta/3 \mbox
{ or } r> \beta/3 \mbox{ and}\quad r\geq(1-\sqrt{1-\beta})^2.
\]
%
%Now if $r\leq\beta/3$, we can take
%$\sqrt{s}=2\sqrt{r}.$ So the best divergence rate of $\Delta_2(s;r,
%under area $r\leq\beta/3$
%in the detectable region
% is of order $L_pp^{\half-\beta+r}.$ If $r>\beta/3,$ the best
%divergence rate of $\Delta_2(s;r,\beta)$ is of order
%$L_pp^{\half-\half(r+\beta)^2/(4r)},$ which is attained at
%$\sqrt{s}=(r+\beta)/(2\sqrt{r}).$

\textit{Case} 4: $s> r$ and $s>(\sqrt{s}-\sqrt{r})^2+\beta$. This
is equivalent to $\sqrt{s}>\max\{(r+\beta)/(2\sqrt{r}),
\sqrt{r}\}$. In this case,
$\mu_{T_{2n},1}(s)-\mu_{T_{2n},0}(s)=L_pp^{1-(\sqrt{s}-\sqrt
{r})^2-\beta}$,\break
$\sigma_{T_{2n},1}(s)=L_pp^{(1-(\sqrt{s}-\sqrt{r})^2-\beta)/2}$. Then
\[
\Delta_2(s;r,\beta)=\frac{\mu_{T_{2n},1}(s)-\mu
_{T_{2n},0}(s)}{\sigma
_{T_{2n},1}(s)}=L_pp^{(1-(\sqrt{s}-\sqrt{r})^2-\beta)/2}.
\]
Hence, it requires that
\[
\sqrt{r}-\sqrt{1-\beta}<\sqrt{s}<\sqrt{r}+\sqrt{1-\beta}.
\]
In order to find an $s$, we
need $\sqrt{r}+\sqrt{1-\beta}>\max\{(r+\beta)/(2\sqrt{r}),
\sqrt{r}\}$. If $\sqrt{r}>(r+\beta)/(2\sqrt{r})$, namely, $r>\beta$,
the above inequality is obviously true. If $r\leq\beta$,
then $\sqrt{r}+\sqrt{1-\beta}>(r+\beta)/(2\sqrt{r})$ is equivalent
to $r>(1-\sqrt{1-\beta})^2$. So the detectable region is
$r>(1-\sqrt{1-\beta})^2$ in this case.
%When $r>\beta,$ we can take $s=r$ such that (\ref{4:eq:4}) is of
%order $L_pp^{(1-\beta)/2}.$ If $(1-\sqrt{1-\beta})^2<r\leq\beta,$
%then the best rate of (\ref{4:eq:4}) is attained at
%$\sqrt{s}=(r+\beta)/(2\sqrt{r}),$ where the best rate is
%$L_pp^{\half-\half(r+\beta)^2/(4r)}.$

In summary of cases 1--4, the union of the detectable regions in the
above four cases is $r>\varrho^*(\beta)$, as illustrated in Figure~\ref{figure5}.
\setcounter{figure}{4}
%
%f5 #&#
\begin{figure}

\includegraphics{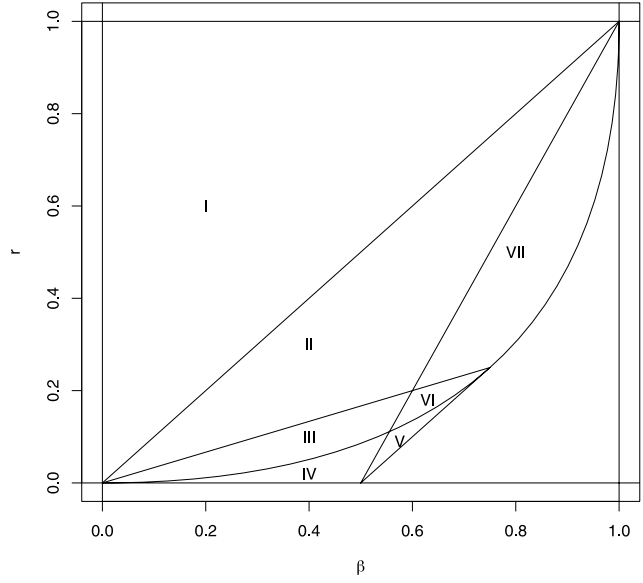}

\caption{The detectable subregions of the $L_2$ threshold test. Case
1: the union of \{I, II, III, IV\}; Case 2, the region is I; Case 3:
the union of \{II, III, IV, V, VI, VII\}; Case 4: the union of \{I, II,
III, VI, VII\}.}

\label{figure5}
\end{figure}

Now we are ready to prove the theorem. We only need to show that the
sum of type I and II
errors of the maximal test goes to 0 when $r>\varrho^*(\beta)$. Because
the maximal test is of
asymptotic $\alpha_n$ level, it suffices to show that
the power goes to 1 in the detectable region
as $n\to\infty$ and $\alpha_n\to0$. Recall that the $\alpha_n$ level
rejection region is $R_{\alpha_n}=\{\hat{\mathcal{M}}_{2n}>\mathcal
{B}_{\alpha_n}\}$.
From Theorem \ref{asy-gumbel}, we notice that $\mathcal{B}_{\alpha_n}=O\{(\log\log
p)^{1/2}\}:=L_p^*$.
%The tightness of $\mathcal{T}_S(s)$ for $s\in\mathcal{S}$
%ensures the existence of $M<\infty$ such that
%$\mathcal{B}_{\alpha_n}<M$ with probability one.
Then, it is sufficient if
%
%eA.6 #&#
\begin{equation}
\label{4:infty} P\bigl(\mathcal{M}_{2n}/L_p^*\to\infty
\bigr)\to1 \qquad\mbox{as } n\to\infty
\end{equation}
at every $(\beta,r)$ in the detectable region. Since $\mathcal
{M}_{2n}\geq\mathcal{T}_{2n}(s)$ for any $s\in\mathcal{S}$,
therefore, (\ref{4:infty}) is true if for any point in the detectable
region, there exists a $\lambda_p(s)=2s\log p$ such that
%
%eA.7 #&#
\begin{equation}
\label{T2ninfinity} \mathcal{T}_{2n}(s)/L_p^*\stackrel{p} {
\to}\infty.
\end{equation}
Therefore, we want to show
%
%eA.8 #&#
\begin{eqnarray}
\label{4:infty:1}&& \frac{T_{2n}(s)-\mu_{T_{2n},0}(s)}{L_p^*\sigma_{Tn,0}(s)}\nonumber\\
&&\qquad= \biggl(\frac
{T_{2n}(s)-\mu_{T_{2n},1}(s)}{L_p^*\sigma_{T_{2n},1}(s)}+\frac{\mu
_{T_{2n},1}(s)-\mu_{T_{2n},0}(s)}{L_p^*\sigma_{T_{2n},1}(s)}
\biggr)\frac
{\sigma_{T_{2n},1}(s)}{\sigma_{T_{2n},0}(s)}
\\
&&\qquad\stackrel{p} {\to} \infty.\nonumber
\end{eqnarray}
%
%For any fixed $M<\infty,$ we want to show
%P(\mathcal{T}(\lambda_p(s))>M)=\Phi\Big(-\frac{\sigma_{Tn,0}}{
%1.
Because ${(T_{2n}(s)-\mu_{T_{2n},1}(s))}/{L_p^*\sigma
_{T_{2n},1}(s)}=o_p(1)$ and $\sigma_{T_{2n},0}(s)\leq\sigma
_{T_{2n},1}(s)$, (\ref{4:infty:1}) is true if
${(\mu_{T_{2n},1}(s)-\mu_{T_{2n},0}(s))}/{L_p^*\sigma
_{T_{2n},1}(s)}\to
\infty$. As we have shown in the early proof, for every $(r,\beta)$ in
the detectable region, there exists an $s$ such that $\frac{\mu
_{T_{2n},1}(s)-\mu_{T_{2n},0}(s)}{L_p\sigma_{T_{2n},1}(s)}\to\infty$
for any slow varying function $L_p$. This
concludes (\ref{T2ninfinity}) and hence (\ref{4:infty}), which
completes the proof of part (i).
%that $T_{2n}(s)\to\infty$ with probability
%one for some $\lambda_p(s)$ and hence the $\mathcal{M}_{2n}$ diverges
%at every
%point in the detectable region, which means that the maximum test can
%attain the optimal detection boundary.

(ii) Note that %$\hat{\mathcal{M}}_{2n}=\mathcal{M}_{2n}\{1+o_p(1)\}$
%and
\[
\hat{\mathcal{M}}_{2n}=\max_{s\in\mathcal{S}_n}
\biggl\{\bigl(\mathcal {T}_{2n,1}(s){R}_2(s)+
\Delta_{2,0}(s;r,\beta)\bigr)\tilde{e}_2(s)+
\frac
{\mu
_{T_{\gamma n},0}(s)-\hat{\mu}_{T_{\gamma n},0}(s)}{\tilde{\sigma
}_{T_{\gamma n}, 0}(s)} \biggr\},
\]
where ${R}_2(s),\tilde{e}_2(s)$ and $\mathcal{T}_{2n,1}(s)$ are defined
in (\ref{th4decom}) and
%
%eA.9 #&#
\begin{eqnarray}
\label{L2-100} \Delta_{2,0}(s;r,\beta)&=&\frac{\mu_{T_{2n},1}(s)-\mu
_{T_{2n},0}(s)}{{\sigma}_{T_{2n}, 0}(s)}\nonumber\\
&=&(s\pi\log
p)^{1/4}(r/s)p^{1/2-\beta+s/2}I(r>s)
\\
&&{} +\frac{s^{1/4}(\pi\log p)^{-1/4}}{2(\sqrt{s}-\sqrt {r})}p^{1/2-(\sqrt{s}-\sqrt{r})^2-\beta+s/2}I(r< s).
\nonumber
\end{eqnarray}

If $r<\varrho^*(\beta)$, then $r<\beta$ and $r<(r+\beta)^2/(4r)$. Hence,
\begin{eqnarray*}
{R}_2(s)=\cases{ %
1 + o(1), &\quad\hspace*{-3pt}$\mbox{if } s\leq r;$
\vspace*{2pt}\cr
1 + o(1), & \quad\hspace*{-3pt}$\mbox{if } r<s \leq
\displaystyle\frac{(r+\beta)^2}{4r};$
\vspace*{2pt}\cr
s^{{1}/{4}}(\sqrt{s}-\sqrt{r})^{-{{1}/{2}}}p^{{{1}/{2}}( 2 \sqrt{s r}-r -
\beta
)}\bigl\{ 1 +
o(1)\bigr\}, & \quad\hspace*{-3pt}$\mbox{if } s>\displaystyle \frac{(r+\beta)^2}{4r}.$}
\end{eqnarray*}
It is also noticed that $r<\varrho^*(\beta)$ implies that $(r+\beta
)^2/(4r)>1$. Therefore, for all $s\in\mathcal{S}_n$, ${R}_2(s)=1+o(1)$.

If $r<\varrho^*(\beta)$, then $r<2\beta-1$. Hence, $\max_{s\leq
r}\Delta
_{2,0}(s;r,\beta)\leq L_pp^{1/2-\beta+r/2}\to0$ as $p(n)\to\infty$.

If $r<\varrho^*(\beta)$ and $r<1/4$, then $r<\beta-1/2$. It follows
that, for all $s>r$,
\[
1/2-(\sqrt{s}-\sqrt{r})^2-\beta+s/2=1/2+r-\beta-\tfrac{1}{2}(
\sqrt {s}-2\sqrt{r})^2\leq1/2+r-\beta<0.
\]
If $r<\varrho^*(\beta)$ and $r>1/4$, then for all $s>r$,
\[
1/2-(\sqrt{s}-\sqrt{r})^2-\beta+s/2%=1/2+r-\beta-\frac{1}{2}(
\leq1/2+r-\beta-\tfrac{1}{2}(1-2\sqrt{r})^2<0.
\]
Hence, $\max_{s> r}\Delta_{2,0}(s;r,\beta)\leq L_pp^{1/2+r-\beta}I\{
r<1/4\}+L_pp^{1-\beta-(1-\sqrt{r})^2}I\{r>1/4\}\to0$ as $p(n)\to
\infty$.
In summary, we have ${R}_2(s)=1+o(1)$ and\break  $\max_{s\in\mathcal
{S}_n}\Delta_{2,0}(s;r,\beta)\to0$ if $r<\varrho^*(\beta)$. Therefore,
together with assumption~(\ref{eq:cri1}), $\hat{\mathcal
{M}}_{2n}=\max_{s\in\mathcal{S}_n}\mathcal{T}_{2n,1}(s)\{1+o_p(1)\}$.

We note that, by employing the same argument of Theorem \ref{asy-gumbel}, it can be
shown that
\[
P \Bigl(a(\log p)\max_{s\in\mathcal{S}}\mathcal{T}_{2n,1}(s)-b(
\log p,\delta)\leq x \Bigr)\to\exp\bigl(-e^{-x}\bigr),
\]
where $\delta$ is defined just above (\ref{Tn1order}). Then the power
of the test
\begin{eqnarray*}
&&P \bigl(\hat{\mathcal{M}}_{2n}>\bigl(\mathcal{E}_{\alpha_n}+b(
\log p,\eta )\bigr)/a(\log p) \bigr)\nonumber
\\
&&\qquad=P \bigl(\hat{\mathcal{M}}_{2n}>\bigl(\mathcal{E}_{\alpha_n}+b(
\log p,\delta )\bigr)/a(\log p) \bigr)\bigl\{1+o(1)\bigr\}\\
&&\qquad=\alpha_n
\bigl\{1+o(1)\bigr\}\to0.\nonumber
\end{eqnarray*}
Thus, the sum of type I and II errors goes to 1. This completes the
proof of part~(ii).
\end{pf*}

\begin{pf*}{Proof of Theorem \ref{th4}}
%To see the power of the maximal $L_2$ test under the alternative
%specified in C.3, we rewrite $\mathcal{M}_{2n}$ as
%Note that
%where $R_2(s)=\frac{ \sigma_{T_{2n} 1}(s)}{\sigma_{T_{2n} 0}(s)}$, $
%&=(s\pi\log p)^{1/4}(r/s)p^{1/2-\beta+s/2}I(r>s)+\frac{s^{1/4}(\pi\log
%p)^{-1/4}}{2(\sqrt{s}-\sqrt{r})}p^{1/2-(\sqrt{s}-\sqrt{r})^2-
We first prove that $\hat{\mathcal{M}}_{\gamma n}\sim\max_{s\in
\mathcal
{S}_n}\Delta_{\gamma,0}(s;r,\beta)$, which will be proved in two parts:
%
%eA.10 #&#
%eA.11 #&#
\begin{eqnarray}
\hat{\mathcal{M}}_{\gamma n}&\sim&{\mathcal{M}}_{\gamma n} \quad\mbox{and}\label{th4part1}
\\
{\mathcal{M}}_{\gamma n}&\sim&\max_{s\in\mathcal{S}_n}\Delta
_{\gamma
,0}(s;r,\beta),\label{th4part2}
\end{eqnarray}
where ${\mathcal{M}}_{\gamma n}=\max_{s\in\mathcal{S}_n}\mathcal
{T}_{\gamma n}(s)=\max_{s\in\mathcal{S}_n} \{\mathcal
{T}_{\gamma
n,1}(s)R_\gamma(s)+\Delta_{\gamma,0}(s;r,\beta) \}$.

To show (\ref{th4part1}), note the decomposition for $\hat{\mathcal
{M}}_{\gamma n}$ in (\ref{th4decom}). Let $\tilde{\mathcal
{M}}_{\gamma
n}=\max_{s\in\mathcal{S}_n}  \{\mathcal{T}_{\gamma n}(s)\tilde
{e}_\gamma(s) \}$. We can first show that $\hat{\mathcal
{M}}_{\gamma
n}\sim\tilde{\mathcal{M}}_{\gamma n}$ because of the following inequality:
\begin{eqnarray*}
&&\tilde{\mathcal{M}}_{\gamma n}- \biggl|\max_{s\in\mathcal{S}_n}
\frac{\mu
_{T_{\gamma n},0}(s)-\hat{\mu}_{T_{\gamma n},0}(s)}{\tilde{\sigma
}_{T_{\gamma n}, 0}(s)} \biggr|
\\
&&\qquad\leq\hat{\mathcal{M}}_{\gamma n} \leq\tilde{
\mathcal{M}}_{\gamma n}+ \biggl|\max_{s\in\mathcal{S}_n} \frac
{\mu_{T_{\gamma n},0}(s)-\hat{\mu}_{T_{\gamma n},0}(s)}{\tilde
{\sigma
}_{T_{\gamma n}, 0}(s)} \biggr|.
\end{eqnarray*}
Under condition (\ref{eq:cri1}), that is, $\max_{s\in\mathcal{S}}\tilde
{\sigma
}^{-1}_{T_{\gamma n}, 0}(s) (\mu_{T_{\gamma n},0}(s)-\hat{\mu
}_{T_{\gamma n},0}(s) )=o(1)$, hence, $\hat{\mathcal{M}}_{\gamma
n}\sim\tilde{\mathcal{M}}_{\gamma n}$.
Second, we can show ${\mathcal{M}}_{\gamma n}\sim\tilde{\mathcal
{M}}_{\gamma n}$. Note the following inequality:
\begin{eqnarray*}
&&\min \Bigl\{{\mathcal{M}}_{\gamma n}\min_{s\in\mathcal{S}_n}\tilde
{e}_\gamma(s),{\mathcal{M}}_{\gamma n}\max_{s\in\mathcal{S}_n}
\tilde {e}_\gamma(s) \Bigr\} \\
&&\qquad\leq\tilde{\mathcal{M}}_{\gamma n}
\leq\max \Bigl\{{\mathcal{M}}_{\gamma n}\min_{s\in\mathcal
{S}_n}
\tilde {e}_\gamma(s),{\mathcal{M}}_{\gamma n}\max
_{s\in\mathcal{S}_n} \tilde {e}_\gamma(s) \Bigr\}.
\end{eqnarray*}
Under conditions (C.1)--(C.4), $\min_{s\in\mathcal{S}_n}\tilde{e}_\gamma
(s)=\max_{s\in\mathcal{S}_n}\tilde{e}_\gamma(s)=1+o(1)$. So we have
\[
\tilde{\mathcal{M}}_{\gamma n}\sim{\mathcal{M}}_{\gamma n}\min
_{s\in
\mathcal{S}_n}\tilde{e}_\gamma(s)\sim{\mathcal{M}}_{\gamma n}
\min_{s\in\mathcal{S}_n}\tilde{e}_\gamma(s)\sim{
\mathcal{M}}_{\gamma n}.
\]
In summary, we have $\hat{\mathcal{M}}_{\gamma n}\sim\tilde
{\mathcal
{M}}_{\gamma n}\sim{\mathcal{M}}_{\gamma n}$. Therefore,
$\hat{\mathcal{M}}_{\gamma n}\sim{\mathcal{M}}_{\gamma n}$.

The path leading to (\ref{th4part2}) is the following. First of all,
it can be shown using an argument similar to the one used in the proof
of Theorem \ref{asy-gumbel} that
\[
P \Bigl(a(\log p)\max_{s\in\mathcal{S}}\mathcal{T}_{\gamma
n,1}(s)-b(
\log p,\delta)\leq x \Bigr)\to\exp\bigl(-e^{-x}\bigr),
\]
where $\delta=\max\{\eta-r+2r\sqrt{1-\eta}-\beta,\eta\}
I(r<1-\eta)+\max
\{1-\beta,\eta\}I(r>1-\eta)$. Thus, for $\gamma=0, 1$ and $2$,
%
%eA.12 #&#
\begin{equation}
\label{Tn1order} \max_{s\in\mathcal{S}}\mathcal{T}_{ \gamma n,1}(s)=O_p
\bigl\{\log ^{1/2}(\log p)\bigr\}.
\end{equation}
%
%=O_p(L_p),$$
%Hence $|\max_{t\in S_n} \mathcal{T}_{\gamma n,1}(t)|=O_p(L_p)$.
%Hence, they grow at a slow varying rate.

Equations (\ref{delta0-sigma-ratio-L2-1}) to (\ref
{delta0-sigma-ratio-L2-8}) in the following reveal that for all $s\in
\mathcal{S}$ and $r>\varrho^*(\beta)$, we can classify $s\in
\mathcal
{S}$ into two sets $\mathcal{S}_1$ and $\mathcal{S}_2$ such that
\begin{eqnarray*}
&&\phantom{i}\mathrm{(i)}\quad \Delta_{\gamma,0}(s;r,\beta)\gg {R}_\gamma(s)\qquad \mbox{for $s\in
\mathcal{S}_1$}
\\
&&\mathrm{(ii)}\quad \Delta_{\gamma,0}(s;r,\beta)\to0\quad \mbox{and}\quad {R}_\gamma
(s)=1+o(1) \qquad\mbox{for $s\in\mathcal{S}_2$},
\end{eqnarray*}
where ``$c\gg d$'' means that $c/d=L_pp^{\xi}$ for some $\xi>0$.
Because $r$ is above the detection boundary $\varrho^*(\beta)$, there
exists at least one $s\in\mathcal{S}_1$ such that $\Delta_{\gamma
,0}(s;r,\beta)\to\infty$.
%This means that the maximum of $\Delta_{\gamma, 0}
Hence,
%
%eA.13 #&#
\begin{equation}
\label{Delta-R-order} \max_{s\in\mathcal{S}}\Delta_{\gamma,0}(s;r,\beta)=
\max_{s\in
\mathcal
{S}_1}\Delta_{\gamma,0}(s;r,\beta)\gg \max
_{s\in\mathcal
{S}}{R}_\gamma(s).
\end{equation}
Namely, the maximum of $\Delta_{\gamma,0}(s;r,\beta)$ is reached on Set
$\mathcal{S}_1$ where $\Delta_{\gamma,0}(s;r,\beta)$ diverges at a much
faster rate than that of $\tilde{R}_\gamma(s)$, if the latter ever diverges.
%As shown in (\ref{max-L2-stat-ineq}) in Appendix
%1-\frac{|\max_{s\in S} \mathcal{T}_{\gamma n,1}(s)R_\gamma(s)|}{\max_{s
%1+\frac{\max_{s\in\mathcal{S}} \mathcal{T}_{\gamma n,1}(s)R_

Let $A(s)=\mathcal{T}_{2n,1}(s){R}_\gamma(s)$. Combining (\ref
{Tn1order}) and (\ref{Delta-R-order}), we have
\[
\Bigl|\max_{s\in\mathcal{S}_n} \mathcal{T}_{\gamma n,1}(s)\Bigr|\Bigl|\max
_{s\in
\mathcal{S}_n}{R}_\gamma(s)\Bigr|=o_p\Bigl\{\max
_{s\in
\mathcal{S}_n}\Delta_{\gamma,0}(s;r,\beta)\Bigr\}.
\]
This implies that $|\max_{s\in\mathcal{S}_n} A(s)|=o_p\{\max_{s\in
\mathcal{S}_n}\Delta_{\gamma,0}(s;r,\beta)\}$.
Together with the following inequality:
\begin{eqnarray*}
\max_{s\in\mathcal{S}_n} \Delta_{\gamma,0}(s;r,\beta)-\Bigl|
\max_{s\in
\mathcal{S}_n} A(s)\Bigr|&\leq&\max_{s\in\mathcal{S}_n}\bigl\{A(s)+
\Delta _{\gamma
,0}(s;r,\beta)\bigr\}
\\
&\leq&\max_{s\in\mathcal{S}_n} \Delta_{\gamma,0}(s;r,\beta)+\max
_{s\in
\mathcal{S}_n} A(s);
\end{eqnarray*}
we conclude that (\ref{th4part2}) holds.

%Next, we are to show the following equivalence in detail
%where `$\sim$' means that both sides of `$\sim$' are at the same order.
%Let $A(s)=\mathcal{T}_{2n,1}(s){R}_\gamma(s)+\tilde{\sigma}^{-1}_{T_{
%n},0}(s)\Big)$. %Because $\Delta_{2,0}(r,t,\beta)+A(t)\geq
%%Since $R_2(t)$ is positive,
%%and $0<\frac{R_2(t)}{\max_s|R_2(s)|}\leq1$, we have
%%$$|\max_{t\in\mathcal{S}_n} \mathcal{T}_{2n,1}(t)\frac{R_2(t)}{
%%which implies that
%%we have $$|\max_{t\in\mathcal{S}_n} A(t)|\leq|\max_t
%Note that

%To establish (\ref{Delta-R-order}), it suffices to attain the following
%&|\max_{s\in\mathcal{S}_n}\tilde{R}_2(s)|=o\{\max_{s\in
%&|\max_{s\in\mathcal{S}_n} \mathcal{T}_{2n,1}(s)|=O_p(L_p).

%By (\ref{delta0-sigma-ratio-L2-1}-\ref{delta0-sigma-ratio-L2-8}),
%$\max_t|R_2(t)|=o\{\max_{t\in\mathcal{S}_n}\Delta_{2,0}(r,t,\beta)\}$
%in cases (i), (ii) and (iii). To show $|\max_{t\in
%we only need to show that $|\max_t \mathcal{T}_{2n,1}(t)|=O_p(L_p)$.

It remains to show the existence of $\mathcal{S}_1$ and $\mathcal{S}_2$
in arriving at (\ref{Delta-R-order}). We only prove it for the $L_2$
test. To complete that, we compare the relative order between $\Delta
_{2,0}(s;r,\beta)$ and ${R}_2(s)$ for three regions above the detection
boundary $\varrho^{\ast}(\beta)$: (i)~$r>\beta$ (ii) $r \in(2
\beta-1,
\beta]$ and (iii) $r\in(\varrho^{\ast}(\beta), 2\beta-1]$. In regions
(i) and (ii) with $r>(1-\sqrt{1-\beta})^2$, we can show that
%
%eA.14 #&#
%eA.15 #&#
\begin{eqnarray}
\label{delta0-sigma-ratio-L2-1}\quad  \Delta_{2,0}(s;r,\beta)&\gg& {R}_2(s)\qquad
\mbox{for $s> 2\beta-1$};
\\
\Delta_{2,0}(s;r,\beta)&\to&0 \quad\mbox{and}\quad {R}_2(s)=1 + o(1)
\qquad\mbox {for $s\leq2\beta-1$}.\label{delta0-sigma-ratio-L2-2}
\end{eqnarray}
%
%where `$a\gg b$' means that $a/b=L_pp^{\alpha}\to\infty$ for a slowing
%varying function $L_p$ and some $\alpha>0$.

In region (ii) with $r<(1-\sqrt{1-\beta})^2$,
%footnote is valid.}
we have
%
%eA.16 #&#
%eA.17 #&#
\begin{eqnarray}\qquad\quad 
\label{delta0-sigma-ratio-L2-3} \Delta_{2,0}(s;r,\beta)&\gg &{R}_2(s)\qquad
\mbox{for $2\beta-1<s\leq (2\sqrt {r}+\sqrt{1+2r-2\beta})^2$},
\\
\Delta_{2,0}(s;r,\beta)&\to&0 \quad\mbox{and}\quad {R}_2(s)=1 + o(1)\qquad
\mbox{for $s\leq2\beta-1$}\label{delta0-sigma-ratio-L2-4}
\nonumber
\\[-8pt]
\\[-8pt]
&& \eqntext{\mbox{and $(2\sqrt{r}+\sqrt{1+2r-2\beta})^2<s< 1$}.}
\end{eqnarray}

For $r\in(\varrho^*(\beta), 2\beta-1]$ in region (iii). If
$r>(1-\sqrt {1-\beta})^2$, define $D_1=(0,(2\sqrt{r}-\sqrt{1+2r-2\beta})^2)$ and
$D_2=((2\sqrt{r}-\sqrt{1+2r-2\beta})^2,1)$. Then it may be shown that
%
%eA.18 #&#
%eA.19 #&#
\begin{eqnarray}
\Delta_{2,0}(s;r,\beta)&\to&0 \quad\mbox{and}\quad {R}_2(s)=1 + o(1)
\qquad\mbox {for $s\in D_1$}; \label{delta0-sigma-ratio-L2-5}
\\
\Delta_{2,0}(s;r,\beta)&\gg& {R}_2(s)\qquad \mbox{for $s\in
D_2$}.\label
{delta0-sigma-ratio-L2-6}
\end{eqnarray}
If $r<(1-\sqrt{1-\beta})^2$, define $D_3=(0,(2\sqrt{r}-\sqrt {1+2r-2\beta
})^2)\cup ((2\sqrt{r}+\break \sqrt{1+2r-2\beta})^2,1)$ and
$D_4=((2\sqrt {r}-\sqrt{1+2r-2\beta})^2, (2\sqrt{r}+ \sqrt{1+2r-2\beta})^2)$.
Then, it
can be shown that
%
%eA.20 #&#
%eA.21 #&#
\begin{eqnarray}
\Delta_{2,0}(s;r,\beta)&\to&0 \quad\mbox{and}\quad {R}_2(s)=1 + o(1)\qquad
\mbox{for $s\in D_3$}; \label{delta0-sigma-ratio-L2-7}
\\
\Delta_{2,0}(s;r,\beta)&\gg& {R}_2(s)\qquad \mbox{for $s\in
D_4$}.\label
{delta0-sigma-ratio-L2-8}
\end{eqnarray}

The results in (\ref{delta0-sigma-ratio-L2-1})--(\ref
{delta0-sigma-ratio-L2-8}) indicate that in each region listed above,
$\max\Delta_{2,0}(s;r,\beta)$ will be attained in situations covered by
(\ref{delta0-sigma-ratio-L2-1}), (\ref{delta0-sigma-ratio-L2-3}),
(\ref
{delta0-sigma-ratio-L2-6}) and (\ref{delta0-sigma-ratio-L2-8}), which
together imply (\ref{Delta-R-order}).

Next, we compute $\Delta_{\gamma,0}(s;r,\beta)$ for the HC ($\gamma=0$)
and the $L_1$ ($\gamma=1$) test. For the HC test, let
$G_{p,1}(s)=P(Y_{i,n}>2s\log p)$. Under assumptions (C.1)--(C.2), %$Y_{i,n}$
%is approximately chi-squared distributed with the noncentral
%component $2r\log p$,
applying the large deviation results [\citet{Petrov}], it may be shown that
%Then the leading order term of $G_{p,1}(s)$ is
\begin{eqnarray*}
G_{p,1}(s)&=&\bigl\{\bigl(2\sqrt{\pi\log p}(\sqrt{s}-\sqrt{r})
\bigr)^{-1}p^{-(\sqrt
{s}-\sqrt{r})^2}\bigr\} \bigl\{1+o(1)\bigr\} \qquad\mbox{if $r<s$
and}
\\
G_{p,1}(s)&=&\bigl\{1-\bigl(2\sqrt{\pi\log p}(\sqrt{r}-\sqrt {s})
\bigr)^{-1}p^{-(\sqrt
{r}-\sqrt{s})^2}\bigr\} \bigl\{1+o(1)\bigr\}\qquad \mbox{if $r>s.$}
\end{eqnarray*}

The mean and variance of $T_{0 n}(s)$ under $H_0$ are
$\mu_{T_{0 n},0}(s)=(\sqrt{s\pi\log p})^{-1}\times p^{1-s}\{1+o(1)\}$ and $
{\sigma}_{T_{0 n},0}^2(s)=(\sqrt{s\pi\log p})^{-1}p^{1-s}\{
1+o(1)\}$ respectively.
The mean and variance of $T_{0 n}(s)$ under the $H_1$ as specified in
(C.4) are, respectively,
\begin{eqnarray*}
\mu_{T_{0 n}, 1} (s)&=& p^{1-\beta
}G_{p,1}(s)+
\bigl(p-p^{1-\beta
}\bigr)2\bar {\Phi}\bigl(\lambda^{1/2}_p(s)
\bigr)\bigl\{1+o(1)\bigr\} \qquad\mbox{and}
\\
% &=&L_pp^{1-s}+p^{1-\beta}I(r>s)+L_pp^{1-(\sqrt{s}-\sqrt{r})^2-
{\sigma}^2_{T_{0 n},1}(s) &=&
p^{1-\beta
}G_{p,1}(s) \bigl(1-G_{p,1}(s)\bigr)\\
&&{}+p
\bigl(1-p^{-\beta}\bigr)2\bar{\Phi}\bigl(\lambda ^{1/2}_p(s)
\bigr) \bigl(1-2\bar{\Phi}\bigl(\lambda^{1/2}_p(s)\bigr)
\bigr).
\end{eqnarray*}
These imply that, up to a factor $\{1+o(1)\}$,
%
%eA.22 #&#
\begin{eqnarray}
&&\mu_{T_{0 n},1}(s)-\mu_{T_{0 n},0}(s)\nonumber\\
&&\qquad=\bigl\{\bigl(2\sqrt{\pi\log p}(
\sqrt {s}-\sqrt{r})\bigr)^{-1}p^{1-\beta-(\sqrt{s}-\sqrt
{r})^2}I(r<s)
\\
&&\hspace*{160pt}\qquad{}+p^{1-\beta
}I(r>s)
\bigr\}\nonumber %\sigma_{\mathcal{Q}_{p,0}}^2(s)&=&(\sqrt{s\pi\log p})^{-1}p^{1-s}
%p})^{-1}p^{1-s}\}\{1+o(1)\}\\
\end{eqnarray}
%
%Let $$\Delta_{HC}(r,s,\beta) = { \mu_{\mathcal{Q}_p,1}(s)-\mu_{
% which is the counterpart of $\Delta_2(s;r,\beta)$ for the
%$L_2$-thresholding formulation.
%The equations given above lead to
%& s\leq(\sqrt{s}-\sqrt{r})^2+\beta; \\
%&s> (\sqrt{s}-\sqrt{r})^2+\beta;\\
% \mbox{and} \\
%& s\leq(\sqrt{s}-\sqrt{r})^2+\beta; \\
%&s> (\sqrt{s}-\sqrt{r})^2+\beta.
and
\[
{R}_0(s)=\cases{ %
1, \qquad
\mbox{if } s\leq(\sqrt{s}-\sqrt{r})^2+\beta;
\vspace*{2pt}\cr
s^{1/4}\bigl|2(\sqrt{s}-\sqrt{r})\bigr|^{-{{1}/{2}}}p^{-{
{1}/{2}}((\sqrt{s}-\sqrt
{r})^2+\beta-s)},\vspace*{2pt}\cr
\quad\hspace*{22pt}\mbox{if } s> (\sqrt{s}-\sqrt{r})^2+\beta.}
\]
Hence,
%
%eA.23 #&#
\begin{eqnarray}
\label{HC-100}\qquad \Delta_{0,0}(s;r,\beta)%:=\frac{\mu_{T_{0n},1}(s)-\mu_{T_{0n},0}(s)}{{
&=&
\frac{s^{1/4}}{2(\sqrt{s}-\sqrt{r})(\pi\log p)^{1/4}}p^{1/2-\beta
-(\sqrt{s}-\sqrt{r})^2+s/2}I(r<s)
\nonumber
\\[-8pt]
\\[-8pt]
\nonumber
&&{} +(s\pi\log p)^{1/4}p^{1/2-\beta+s/2}I(r>s).
\end{eqnarray}

%%%%%%
%% Mean and variance for L-1 test statistic
%%%%%%

For the $L_1$ test, the mean and variances of $T_{1n}(s)$ under $H_1$
specified in (C.4) are, respectively, up to a factor $1+o(1)$,
\begin{eqnarray*}
\mu_{T_{1n},1}(s)&=&\frac{\sqrt{s}}{\sqrt{2\pi}(\sqrt{s}-\sqrt {r})}p^{1-\beta-(\sqrt{s}-\sqrt{r})^2}I(r<s)
\\
&&{} +(\sqrt{2r\log p}) p^{1-\beta}I(r>s)+\sqrt{2/\pi}p^{1-s}%\{1+o(1)
\qquad\mbox{and}
\\
{\sigma}^2_{T_{1n},1}(s)&=&\frac{s\sqrt{\log p}}{\sqrt{\pi}(\sqrt {s}-\sqrt{r})}p^{1-\beta-(\sqrt{s}-\sqrt{r})^2}I(r<s)
+p^{1-\beta}I(r>s)\\
&&{}+2\sqrt{(s/\pi)\log p}p^{1-s}. %\{1+o(1)\}.
\end{eqnarray*}
It follows that, up to a factor $1+o(1)$,
%
%eA.24 #&#
\begin{eqnarray}
\mu_{T_{1n},1}(s)-\mu_{T_{1n},0}(s)& =&\frac{\sqrt{s}}{\sqrt{2\pi}(\sqrt{s}-\sqrt{r})}p^{1-\beta
-(\sqrt
{s}-\sqrt{r})^2}I(r<s)
\nonumber
\\[-8pt]
\\[-8pt]
\nonumber
&&{}+(
\sqrt{2r\log p}) p^{1-\beta}I(r>s) %\{1+o(1)\}
\end{eqnarray}
%
% \mbox{and} \\
%& s\leq(\sqrt{s}-\sqrt{r})^2+\beta; \\
%&s> (\sqrt{s}-\sqrt{r})^2+\beta.
and
\[
{R}_1(s)=\cases{ %
1, \qquad
\mbox{if } s\leq r \mbox{ and } s\leq\beta;
\vspace*{2pt}\cr
\displaystyle(\sqrt{2})^{-1}\biggl(\frac{s}{\pi}\biggr)^{-{1}/{4}}(\log
p)^{-{1}/{4}}p^{(s-\beta
)/2}, %\{1+o(1)\}
\vspace*{2pt}\cr
\qquad\quad \mbox{if }  s\leq r \mbox{ and } s
\geq\beta;
\vspace*{2pt}\cr
1, %1+o(1)
 \qquad\mbox{if } s> r \mbox{ and } s\leq(\sqrt{s}-\sqrt{r})^2+
\beta;
\vspace*{2pt}\cr
s^{{1}/{4}}(2\sqrt{s}-2\sqrt{r})^{-{{1}/{2}}}p^{-{
{1}/{2}}((\sqrt
{s}-\sqrt
{r})^2+\beta-s)},\vspace*{2pt}\cr
 \qquad\quad\mbox{if } s> r \mbox{ and } s> (\sqrt{s}-\sqrt{r})^2+
\beta.}
\]
%
%where $R_1(s)=\frac{\sigma_{T_{1n},1}(s)}{\sigma_{T_{1n},0}(s)}$.
Therefore,
%
%eA.25 #&#
\begin{eqnarray}\label{L1-100}
\qquad\Delta_{1,0}(s;r,\beta)%:=\frac{\mu_{T_{1n},1}(s)-\mu_{T_{1n},0}(s)}{{
&=&
\frac{s^{{1}/{4}}}{2(\pi\log p)^{{1}/{4}}(\sqrt {s}-\sqrt {r})}p^{1/2-\beta-(\sqrt{s}-\sqrt{r})^2+s/2}I(r<s)
\nonumber
\\[-8pt]
\\[-8pt]
\nonumber
&&{} +(s\pi\log p)^{{1}/{4}}(r/s)^{{1}/{4}} p^{1/2-\beta
+s/2}I(r>s).
\end{eqnarray}

Replicating the above proof for the $L_2$ test, it can be shown that,
for $\gamma=0$ and~1,
\[
\hat{\mathcal{M}}_{\gamma n}\sim\max_{s\in\mathcal{S}_n}
\Delta _{\gamma
,0}(s;r,\beta).
\]
%
%where
%&=\frac{s^{\fourth}}{2(\pi\log p)^{\fourth}(\sqrt{s}-\sqrt{r})}p^{1/2-
%& +(s\pi\log p)^{\fourth}(r/s)^{\fourth} p^{1/2-\beta+s/2}I(r>s)
%%\end{eqnarray}
%%A very similar analysis to above derivation shows that
%%\begin{eqnarray}
%%\label{HC-equi}
%%\mathcal{M}_{0n}\sim\max_{s\in\mathcal{S}_n}\Delta_{0,0}(s;r,\beta)
%%\end{eqnarray}
%%where
%%\begin{eqnarray}
%p)^{1/4}}p^{1/2-\beta-(\sqrt{s}-\sqrt{r})^2+s/2}I(r<s)\nn\\
%& +(s\pi\log p)^{1/4}p^{1/2-\beta+s/2}I(r>s).

At last, we will compare $\max_{s\in\mathcal{S}_n}\Delta_{\gamma
,0}(s;r,\beta)$ for $\gamma=0,1$ and 2 when $r>2\beta-1$.
Let $s^*_n=\arg\max\{s\dvtx s\in\mathcal{S}_n\cap(2\beta-1,r)\}$ be a
threshold in $(2\beta-1,r)$ that is closest to $r$. Then the maximal
value of $\Delta_{\gamma,0}(s,r,\beta)$ over $\mathcal{S}_n$ is
attained at $s^*_n$. Note that such $s_n^*$ exists with probability 1.
To show this point, it is enough to show that $\mathcal{S}_n\cap
(2\beta
-1,r)\neq\varnothing$, which is equivalent to showing that $P(\bigcup_{i=1}^p\{Y_{i,n}\in((4\beta-2)\log p, 2r\log p)\})\to1$. Let $\{
k_1,\ldots,k_q\}\in(1,\ldots,p)$ be a sub-sequence such that $q\to
\infty$ and $k_{\min}=\min_j|k_j-k_{j-1}|\to\infty$. Let
$D_n=\prod_{i=k_1}^{k_q}P(\{Y_{i,n}\in((4\beta-2)\log p, 2r\log p)^c\})-P(\bigcap_{i=k_1}^{k_q}\{Y_{i,n}\in((4\beta-2)\log p, 2r\log p)^c\})$. By
mixing assumption (C.5) and the triangle inequality, it can be seen
that $|D_n|\leq q\alpha_Z(k_{\min})\to0$ as $n\to\infty$. Then it
follows that
\begin{eqnarray*}
&&P\Biggl(\bigcup_{i=1}^p\bigl
\{Y_{i,n}\in\bigl((4\beta-2)\log p, 2r\log p\bigr)\bigr\}\Biggr)
\\
&&\qquad\geq P\Biggl(\bigcup_{i=k_1}^{k_{q}}\bigl
\{Y_{i,n}\in\bigl((4\beta-2)\log p, 2r\log p\bigr)\bigr\} \Biggr)
\\
&&\qquad=1-P\Biggl(\bigcap_{i=k_1}^{k_q}\bigl
\{Y_{i,n}\in\bigl((4\beta-2)\log p, 2r\log p\bigr)^c\bigr\}
\Biggr)
\\
&&\qquad=1-\prod_{i=k_1}^{k_q}P\bigl(\bigl
\{Y_{i,n}\in\bigl((4\beta-2)\log p, 2r\log p\bigr)^c\bigr\}
\bigr)+D_n\to1,
\end{eqnarray*}
where we used $P(\{Y_{i,n}\in((4\beta-2)\log p, 2r\log p)^c\})<1$ for
all $i=1,\ldots,p$. Comparing (\ref{L2-100}), (\ref{HC-100}) and
(\ref
{L1-100}), we see that $\Delta_{0,0}(s_n^*;r,\beta)< \Delta
_{1,0}(s_n^*;\break r,\beta)<\Delta_{2,0}(s_n^*;r,\beta)$.

It follows that, for $r>2\beta-1$,
\[
\max_{s\in\mathcal{S}_n}\Delta_{0,0}(s;r,\beta)<
\max_{s\in
\mathcal
{S}_n}\Delta_{1,0}(s;r,\beta) < \max
_{s\in\mathcal{S}_n}\Delta _{2,0}(s;r,\beta).
\]
Therefore, asymptotically with probability 1,
$\hat{\mathcal{M}}_{0 n} < \hat{\mathcal{M}}_{1n}< \hat{\mathcal{M}}_{2n}$,
which results in $\Omega_{0}(r,\beta)\leq\Omega_1(r,\beta)\leq
\Omega
_2(r,\beta)$. This completes the proof.
\end{pf*}
\end{appendix}

\section*{Acknowledgments}
The authors thank
the Editor, an Associate Editor and two referees for insightful and
constructive comments which have improved the presentation of the paper.
We are also very grateful to Dr. Jiashun Jin and Dr. Cun-Hui
Zhang for stimulating discussions.

\begin{supplement}[id=suppA]
\stitle{A supplement to ``Tests alternative to
higher criticism for high-dimensional means under sparsity and
column-wise dependence''\\}
\slink[doi]{10.1214/13-AOS1168SUPP} %[doi,text={...}] - jei reikia
%suskaldyti doi
\sdatatype{.pdf}
\sfilename{aos1168\_supp.pdf}
\sdescription{The supplementary material contains proofs for
Proposition~\ref{chap4-cor2} and Theorem \ref{th1} in Section~\ref{sec2}.}
\end{supplement}

% imsref loaded by akundreckaite, 2013-11-13 12:12:34
%
% imsref loaded by akundreckaite, 2013-11-13 13:11:27

%


\printaddresses


\begin{thebibliography}{31}
% BibTex style file: ims.bst, 2013-01-28
% Default style options (sort=0,type=number).
% Used options (sort=1,type=nameyear).

%b1 #&#
\bibitem[\protect\citeauthoryear{Andrews and Pollard}{1994}]{Andrews}
\begin{barticle}[auto:STB|2013/10/14|10:36:11]
\bauthor{\bsnm{Andrews},~\bfnm{D.}\binits{D.}} \AND
\bauthor{\bsnm{Pollard},~\bfnm{D.}\binits{D.}}
(\byear{1994}).
\btitle{An introduction to functional central limit theorems for dependent
stochastic processes}.
\bjournal{Int. Statist. Rev.}
\bvolume{62}
\bpages{119--132}.
\bptok{imsref}%
\end{barticle}
\endbibitem

%b2 #&#
\bibitem[\protect\citeauthoryear{Arias-Castro, Bubeck and
Lugosi}{2012a}]{Arias-Castro2012a}
\begin{barticle}[mr]
\bauthor{\bsnm{Arias-Castro},~\bfnm{Ery}\binits{E.}},
\bauthor{\bsnm{Bubeck},~\bfnm{S{\'e}bastien}\binits{S.}} \AND
\bauthor{\bsnm{Lugosi},~\bfnm{G{\'a}bor}\binits{G.}}
(\byear{2012}a).
\btitle{Detection of correlations}.
\bjournal{Ann. Statist.}
\bvolume{40}
\bpages{412--435}.
\bid{doi={10.1214/11-AOS964}, issn={0090-5364}, mr={3014312}}
\bptok{imsref}%
\end{barticle}
\endbibitem

%b3 #&#
\bibitem[\protect\citeauthoryear{Arias-Castro, Bubeck and
Lugosi}{2012b}]{Arias-Castro2012b}
\begin{bmisc}[auto:STB|2013/10/14|10:36:11]
\bauthor{\bsnm{Arias-Castro},~\bfnm{E.}\binits{E.}},
\bauthor{\bsnm{Bubeck},~\bfnm{S.}\binits{S.}} \AND
\bauthor{\bsnm{Lugosi},~\bfnm{G.}\binits{G.}}
(\byear{2012}b).
\bhowpublished{Detecting positive correlations in a multivariate sample.
Available at \arxivurl{arXiv:1202.5536v1} [math.ST]}.
\bptok{imsref}%
\end{bmisc}
\endbibitem

%b4 #&#
\bibitem[\protect\citeauthoryear{Bai and Saranadasa}{1996}]{BaiSara}
\begin{barticle}[mr]
\bauthor{\bsnm{Bai},~\bfnm{Zhidong}\binits{Z.}} \AND
\bauthor{\bsnm{Saranadasa},~\bfnm{Hewa}\binits{H.}}
(\byear{1996}).
\btitle{Effect of high dimension: By an example of a two sample problem}.
\bjournal{Statist. Sinica}
\bvolume{6}
\bpages{311--329}.
\bid{issn={1017-0405}, mr={1399305}}
\bptok{imsref}%
\end{barticle}
\endbibitem

%b5 #&#
\bibitem[\protect\citeauthoryear{Bradley}{2005}]{Bradley}
\begin{barticle}[mr]
\bauthor{\bsnm{Bradley},~\bfnm{Richard~C.}\binits{R.~C.}}
(\byear{2005}).
\btitle{Basic properties of strong mixing conditions. {A} survey and some open
questions}.
\bjournal{Probab. Surv.}
\bvolume{2}
\bpages{107--144}.
\bid{doi={10.1214/154957805100000104}, issn={1549-5787}, mr={2178042}}
\bptok{imsref}%
\end{barticle}
\endbibitem

%b6 #&#
\bibitem[\protect\citeauthoryear{Cai, Jeng and Jin}{2011}]{Cai2011}
\begin{barticle}[mr]
\bauthor{\bsnm{Cai},~\bfnm{T.~Tony}\binits{T.~T.}},
\bauthor{\bsnm{Jeng},~\bfnm{X.~Jessie}\binits{X.~J.}} \AND
\bauthor{\bsnm{Jin},~\bfnm{Jiashun}\binits{J.}}
(\byear{2011}).
\btitle{Optimal detection of heterogeneous and heteroscedastic mixtures}.
\bjournal{J. R. Stat. Soc. Ser. B Stat. Methodol.}
\bvolume{73}
\bpages{629--662}.
\bid{doi={10.1111/j.1467-9868.2011.00778.x}, issn={1369-7412}, mr={2867452}}
\bptok{imsref}%
\end{barticle}
\endbibitem

%b7 #&#
\bibitem[\protect\citeauthoryear{Cai and Wu}{2012}]{Cai2012}
\begin{bmisc}[auto:STB|2013/10/14|10:36:11]
\bauthor{\bsnm{Cai},~\bfnm{T.}\binits{T.}} \AND
\bauthor{\bsnm{Wu},~\bfnm{Y.}\binits{Y.}}
(\byear{2012}).
\bhowpublished{Optimal detection for sparse mixtures. Unpublished manuscript.}
\bptok{imsref}%
\end{bmisc}
\endbibitem

%b8 #&#
\bibitem[\protect\citeauthoryear{Chen and Qin}{2010}]{ChenQin}
\begin{barticle}[mr]
\bauthor{\bsnm{Chen},~\bfnm{Song~Xi}\binits{S.~X.}} \AND
\bauthor{\bsnm{Qin},~\bfnm{Ying-Li}\binits{Y.-L.}}
(\byear{2010}).
\btitle{A two-sample test for high-dimensional data with applications to
gene-set testing}.
\bjournal{Ann. Statist.}
\bvolume{38}
\bpages{808--835}.
\bid{doi={10.1214/09-AOS716}, issn={0090-5364}, mr={2604697}}
\bptok{imsref}%
\end{barticle}
\endbibitem

%b9 #&#
\bibitem[\protect\citeauthoryear{Delaigle and Hall}{2009}]{Delaigle2009}
\begin{bincollection}[mr]
\bauthor{\bsnm{Delaigle},~\bfnm{Aurore}\binits{A.}} \AND
\bauthor{\bsnm{Hall},~\bfnm{Peter}\binits{P.}}
(\byear{2009}).
\btitle{Higher criticism in the context of unknown distribution,
no-nindependence and classification}.
In \bbooktitle{Perspectives in Mathematical Sciences. {I}}.
\bseries{Stat. Sci. Interdiscip. Res.}
\bvolume{7}
\bpages{109--138}.
\bpublisher{World Sci. Publ.}, \blocation{Hackensack, NJ}.
\bid{doi={10.1142/9789814273633_0006}, mr={2581742}}
\bptok{imsref}%
\end{bincollection}
\endbibitem

%b10 #&#
\bibitem[\protect\citeauthoryear{Delaigle, Hall and Jin}{2011}]{Delaigle}
\begin{barticle}[mr]
\bauthor{\bsnm{Delaigle},~\bfnm{Aurore}\binits{A.}},
\bauthor{\bsnm{Hall},~\bfnm{Peter}\binits{P.}} \AND
\bauthor{\bsnm{Jin},~\bfnm{Jiashun}\binits{J.}}
(\byear{2011}).
\btitle{Robustness and accuracy of methods for high dimensional data analysis
based on {S}tudent's {$t$}-statistic}.
\bjournal{J. R. Stat. Soc. Ser. B Stat. Methodol.}
\bvolume{73}
\bpages{283--301}.
\bid{doi={10.1111/j.1467-9868.2010.00761.x}, issn={1369-7412}, mr={2815777}}
\bptok{imsref}%
\end{barticle}
\endbibitem

%b11 #&#
\bibitem[\protect\citeauthoryear{Donoho and Jin}{2004}]{DonohoJin}
\begin{barticle}[mr]
\bauthor{\bsnm{Donoho},~\bfnm{David}\binits{D.}} \AND
\bauthor{\bsnm{Jin},~\bfnm{Jiashun}\binits{J.}}
(\byear{2004}).
\btitle{Higher criticism for detecting sparse heterogeneous mixtures}.
\bjournal{Ann. Statist.}
\bvolume{32}
\bpages{962--994}.
\bid{doi={10.1214/009053604000000265}, issn={0090-5364}, mr={2065195}}
\bptok{imsref}%
\end{barticle}
\endbibitem

%b12 #&#
\bibitem[\protect\citeauthoryear{Donoho and Jin}{2008}]{DonohoJin08}
\begin{barticle}[pbm]
\bauthor{\bsnm{Donoho},~\bfnm{David}\binits{D.}} \AND
\bauthor{\bsnm{Jin},~\bfnm{Jiashun}\binits{J.}}
(\byear{2008}).
\btitle{Higher criticism thresholding: Optimal feature selection when useful
features are rare and weak}.
\bjournal{Proc. Natl. Acad. Sci. USA}
\bvolume{105}
\bpages{14790--14795}.
\bid{doi={10.1073/pnas.0807471105}, issn={1091-6490}, pii={0807471105},
pmcid={2553037}, pmid={18815365}}
\bptok{imsref}%
\end{barticle}
\endbibitem

%b13 #&#
\bibitem[\protect\citeauthoryear{Donoho and Johnstone}{1994}]{DonohoJohnstone}
\begin{barticle}[mr]
\bauthor{\bsnm{Donoho},~\bfnm{David~L.}\binits{D.~L.}} \AND
\bauthor{\bsnm{Johnstone},~\bfnm{Iain~M.}\binits{I.~M.}}
(\byear{1994}).
\btitle{Ideal spatial adaptation by wavelet shrinkage}.
\bjournal{Biometrika}
\bvolume{81}
\bpages{425--455}.
\bid{doi={10.1093/biomet/81.3.425}, issn={0006-3444}, mr={1311089}}
\bptok{imsref}%
\end{barticle}
\endbibitem

%b14 #&#
\bibitem[\protect\citeauthoryear{Doukhan}{1994}]{Doukhan}
\begin{bbook}[mr]
\bauthor{\bsnm{Doukhan},~\bfnm{Paul}\binits{P.}}
(\byear{1994}).
\btitle{Mixing: Properties and Examples}.
\bseries{Lecture Notes in Statistics}
\bvolume{85}.
\bpublisher{Springer}, \blocation{New York}.
\bid{doi={10.1007/978-1-4612-2642-0}, mr={1312160}}
\bptok{imsref}%
\end{bbook}
\endbibitem

%b15 #&#
\bibitem[\protect\citeauthoryear{Fan}{1996}]{Fan1996}
\begin{barticle}[mr]
\bauthor{\bsnm{Fan},~\bfnm{Jianqing}\binits{J.}}
(\byear{1996}).
\btitle{Test of significance based on wavelet thresholding and {N}eyman's
truncation}.
\bjournal{J. Amer. Statist. Assoc.}
\bvolume{91}
\bpages{674--688}.
\bid{doi={10.2307/2291663}, issn={0162-1459}, mr={1395735}}
\bptok{imsref}%
\end{barticle}
\endbibitem

%b16 #&#
\bibitem[\protect\citeauthoryear{Hall and Jin}{2008}]{HallJin}
\begin{barticle}[mr]
\bauthor{\bsnm{Hall},~\bfnm{Peter}\binits{P.}} \AND
\bauthor{\bsnm{Jin},~\bfnm{Jiashun}\binits{J.}}
(\byear{2008}).
\btitle{Properties of higher criticism under strong dependence}.
\bjournal{Ann. Statist.}
\bvolume{36}
\bpages{381--402}.
\bid{doi={10.1214/009053607000000767}, issn={0090-5364}, mr={2387976}}
\bptok{imsref}%
\end{barticle}
\endbibitem

%b17 #&#
\bibitem[\protect\citeauthoryear{Hall and Jin}{2010}]{HallJin2010}
\begin{barticle}[mr]
\bauthor{\bsnm{Hall},~\bfnm{Peter}\binits{P.}} \AND
\bauthor{\bsnm{Jin},~\bfnm{Jiashun}\binits{J.}}
(\byear{2010}).
\btitle{Innovated higher criticism for detecting sparse signals in correlated
noise}.
\bjournal{Ann. Statist.}
\bvolume{38}
\bpages{1686--1732}.
\bid{doi={10.1214/09-AOS764}, issn={0090-5364}, mr={2662357}}
\bptok{imsref}%
\end{barticle}
\endbibitem

%b18 #&#
\bibitem[\protect\citeauthoryear{Ingster}{1997}]{Ingster}
\begin{barticle}[mr]
\bauthor{\bsnm{Ingster},~\bfnm{Yu.~I.}\binits{Y.~I.}}
(\byear{1997}).
\btitle{Some problems of hypothesis testing leading to infinitely divisible
distributions}.
\bjournal{Math. Methods Statist.}
\bvolume{6}
\bpages{47--69}.
\bid{issn={1066-5307}, mr={1456646}}
\bptok{imsref}%
\end{barticle}
\endbibitem

%b19 #&#
\bibitem[\protect\citeauthoryear{Jing, Shao and Zhou}{2008}]{Jing}
\begin{barticle}[mr]
\bauthor{\bsnm{Jing},~\bfnm{Bing-Yi}\binits{B.-Y.}},
\bauthor{\bsnm{Shao},~\bfnm{Qi-Man}\binits{Q.-M.}} \AND
\bauthor{\bsnm{Zhou},~\bfnm{Wang}\binits{W.}}
(\byear{2008}).
\btitle{Towards a universal self-normalized moderate deviation}.
\bjournal{Trans. Amer. Math. Soc.}
\bvolume{360}
\bpages{4263--4285}.
\bid{doi={10.1090/S0002-9947-08-04402-4}, issn={0002-9947}, mr={2395172}}
\bptok{imsref}%
\end{barticle}
\endbibitem

%b20 #&#
\bibitem[\protect\citeauthoryear{Joe}{1997}]{Joe}
\begin{bbook}[mr]
\bauthor{\bsnm{Joe},~\bfnm{Harry}\binits{H.}}
(\byear{1997}).
\btitle{Multivariate Models and Dependence Concepts}.
\bseries{Monographs on Statistics and Applied Probability}
\bvolume{73}.
\bpublisher{Chapman \& Hall}, \blocation{London}.
\bid{mr={1462613}}
\bptok{imsref}%
\end{bbook}
\endbibitem

%b21 #&#
\bibitem[\protect\citeauthoryear{Leadbetter, Lindgren and
Rootz{\'e}n}{1983}]{Leadbetter}
\begin{bbook}[mr]
\bauthor{\bsnm{Leadbetter},~\bfnm{M.~R.}\binits{M.~R.}},
\bauthor{\bsnm{Lindgren},~\bfnm{Georg}\binits{G.}} \AND
\bauthor{\bsnm{Rootz{\'e}n},~\bfnm{Holger}\binits{H.}}
(\byear{1983}).
\btitle{Extremes and Related Properties of Random Sequences and Processes}.
\bpublisher{Springer}, \blocation{New York}.
\bid{mr={0691492}}
\bptok{imsref}%
\end{bbook}
\endbibitem

%b22 #&#
\bibitem[\protect\citeauthoryear{Petrov}{1995}]{Petrov}
\begin{bbook}[mr]
\bauthor{\bsnm{Petrov},~\bfnm{Valentin~V.}\binits{V.~V.}}
(\byear{1995}).
\btitle{Limit Theorems of Probability Theory: Sequences of Independent Random
Variables}.
\bseries{Oxford Studies in Probability}
\bvolume{4}.
\bpublisher{Oxford Univ. Press}, \blocation{New York}.
\bid{mr={1353441}}
\bptok{imsref}%
\end{bbook}
\endbibitem

%b23 #&#
\bibitem[\protect\citeauthoryear{Pisier}{1983}]{Pisier}
\begin{bincollection}[mr]
\bauthor{\bsnm{Pisier},~\bfnm{Gilles}\binits{G.}}
(\byear{1983}).
\btitle{Some applications of the metric entropy condition to harmonic
analysis}.
In \bbooktitle{Banach Spaces, Harmonic Analysis, and Probability Theory
({S}torrs, {C}onn., 1980/1981)}.
\bseries{Lecture Notes in Math.}
\bvolume{995}
\bpages{123--154}.
\bpublisher{Springer}, \blocation{Berlin}.
\bid{doi={10.1007/BFb0061891}, mr={0717231}}
\bptok{imsref}%
\end{bincollection}
\endbibitem

%b24 #&#
\bibitem[\protect\citeauthoryear{Shao}{1997}]{Shao}
\begin{barticle}[mr]
\bauthor{\bsnm{Shao},~\bfnm{Qi-Man}\binits{Q.-M.}}
(\byear{1997}).
\btitle{Self-normalized large deviations}.
\bjournal{Ann. Probab.}
\bvolume{25}
\bpages{285--328}.
\bid{doi={10.1214/aop/1024404289}, issn={0091-1798}, mr={1428510}}
\bptok{imsref}%
\end{barticle}
\endbibitem

%%b25 #&#
%(\byear{1986}).

%b26 #&#
\bibitem[\protect\citeauthoryear{Sibuya}{1960}]{Sibuya}
\begin{barticle}[mr]
\bauthor{\bsnm{Sibuya},~\bfnm{Masaaki}\binits{M.}}
(\byear{1960}).
\btitle{Bivariate extreme statistics. {I}}.
\bjournal{Ann. Inst. Statist. Math. Tokyo}
\bvolume{11}
\bpages{195--210}.
\bid{issn={0020-3157}, mr={0115241}}
\bptok{imsref}%
\end{barticle}
\endbibitem

%b27 #&#
\bibitem[\protect\citeauthoryear{Tukey}{1976}]{Tukey}
\begin{bmisc}[auto:STB|2013/10/14|10:36:11]
\bauthor{\bsnm{Tukey},~\bfnm{J.~W.}\binits{J.~W.}}
(\byear{1976}).
\bhowpublished{T13 N: The higher criticism. Course Notes, Statistics 411,
Princeton Univ}.
\bptok{imsref}%
\end{bmisc}
\endbibitem

%b28 #&#
\bibitem[\protect\citeauthoryear{van~der Vaart and Wellner}{1996}]{Wellner}
\begin{bbook}[mr]
\bauthor{\bparticle{van~der} \bsnm{Vaart},~\bfnm{Aad~W.}\binits{A.~W.}} \AND
\bauthor{\bsnm{Wellner},~\bfnm{Jon~A.}\binits{J.~A.}}
(\byear{1996}).
\btitle{Weak Convergence and Empirical Processes: With Applications to
Statistics}.
\bpublisher{Springer}, \blocation{New York}.
\bid{mr={1385671}}
\bptok{imsref}%
\end{bbook}
\endbibitem

%b29 #&#
\bibitem[\protect\citeauthoryear{Wang and Hall}{2009}]{WangHall}
\begin{barticle}[mr]
\bauthor{\bsnm{Wang},~\bfnm{Qiying}\binits{Q.}} \AND
\bauthor{\bsnm{Hall},~\bfnm{Peter}\binits{P.}}
(\byear{2009}).
\btitle{Relative errors in central limit theorems for {S}tudent's {$t$}
statistic, with applications}.
\bjournal{Statist. Sinica}
\bvolume{19}
\bpages{343--354}.
\bid{issn={1017-0405}, mr={2487894}}
\bptok{imsref}%
\end{barticle}
\endbibitem

%b30 #&#
\bibitem[\protect\citeauthoryear{Zhong, Chen and Xu}{2013}]{ZhenChenXu}
\begin{bmisc}[auto:STB|2013/10/14|10:36:11]
\bauthor{\bsnm{Zhong},~\bfnm{P.~S.}\binits{P.~S.}},
\bauthor{\bsnm{Chen},~\bfnm{S.~X.}\binits{S.~X.}} \AND
\bauthor{\bsnm{Xu},~\bfnm{M.}\binits{M.}}
(\byear{2013}).
\bhowpublished{Supplement to ``Tests alternative to higher criticism for high
dimensional means under sparsity and column-wise dependence.''
DOI:\doiurl{10.1214/13-AOS1168SUPP}.}
\bptok{imsref}%
\end{bmisc}
\endbibitem

\end{thebibliography}
\end{document}